# WEYL-HEISENBERG FRAMES, TRANSLATION INVARIANT SYSTEMS AND THE WALNUT REPRESENTATION

Peter G. Casazza, Ole Christensen and A.J.E.M. Janssen

ABSTRACT. We present a comprehensive analysis of the convergence properties of the frame operators of Weyl-Heisenberg systems and shift-invariant systems, and relate these to the convergence of the Walnut representation. We give a deep analysis of necessary conditions and sufficient conditions for convergence of the frame operator. We show that symmetric, norm and unconditional convergence of the Walnut series are all different, but that weak and norm convergence are the same, while there are WH-systems for which the Walnut representation has none of these convergence properties. We make a detailed study of the CC-condition (a sufficient condition for WH-systems to have finite upper frame bounds), and show that (for ab rational) a uniform version of this passes to the Wexler-Raz dual. We also show that a condition of Tolimieri and Orr implies the uniform CC-condition. We obtain stronger results in the case when $(g, a, b)$ is a WH-system and ab is rational. For example, if ab is rational, then the CC-condition becomes equivalent to the unconditional convergence of the Walnut representation - even in a more general setting. Many of the results are generalized to shift-invariant systems. We give classifications for numerous important classes of WH-systems including: (1) The WH-systems for which the frame operator extends to a bounded operator on $L^p(R)$, for all $1 \leq p \leq \infty$; (2) The WH-systems for which the frame operator extends to a bounded operator on the Wiener amalgam space; (3) The families of frames which have the same frame operator.

## 1. Introduction

Well after introduction by D. Gabor [12] in 1946, Weyl-Heisenberg systems (also known as Gabor systems, see the Introduction in [10] for the history) nowadays play a central role in e.g. signal processing, image and data compression. However, major difficulties can arise when one tries to compute the frame operator for these systems. Walnut [24] gave a series representation for the frame operator under the assumption that the generating function for the frame has rapid decay. In this case, the Walnut series converges rapidly in norm. The Walnut representation of the frame operator can often be much easier to compute and work with. But the

The first author was supported by NSF DMS 970618.

Typeset by $\mathcal{A}_{\mathcal{M}}\mathcal{S}$-TeX

rapid decay assumptions limit its use. Recently, Casazza and Christensen [2] gave much weaker assumptions on the window function $g$ which still ensures that the corresponding Weyl-Heisenberg system has a finite upper frame bound. As we will see, this condition also has important implications for the Walnut representation of the frame operator.

Until now, there has not been a detailed study of the frame operator for WH-systems and its relationship to the Walnut representation, especially, in the general case where we do not put rapid decay assumptions on $g$. This paper is a comprehensive study of the convergence properties of the frame operator for WH-systems and shift-invariant systems, as well as the Walnut representation, and the relationships between them. The body of the paper is divided into seven sections. Section 2 contains the basic tools we need for this study including the basic definitions, the WH-frame identity, and a detailed outline of the types of convergence (weak, norm, and unconditional) which we work with in the paper. Section 3 introduces the Zak transform, its basic properties. Then, for the case $a = b = 1$, we use the Zak transform to produce necessary conditions for the existence of finite upper frame bounds. This section also contains the basic Zak transform construction for examples which will be used throughout the paper. Section 4 is a detailed study of sufficient conditions for the convergence of the frame operator. We consider a much weaker condition than previously used and show that this condition is maintained under application of the frame operator to the window function even for general shift-invariant systems.. With slightly stronger assumptions we show that this condition can be passed to the Wexler-Raz dual window. In Section 5 we analyze symmetric weak and norm convergence of the Walnut representation of the frame operator. We give classifications for the families of WH-systems for which the Walnut representation has each of these convergence properties. One consequence is that weak and norm convergence of the Walnut representation are equivalent. We also give examples to show that the other forms of convergence are all different. In Section 6 we consider unconditional convergence of the Walnut series. In this case, we give several classifications of this property and show that the conditions considered in Section 4 imply unconditional convergence of the Walnut series. We also produce a large class of functions for which the Walnut series converges with respect to every WH-system. In Section 7 we classify the WH-systems which have the property that their frame operator extends to be a bounded linear operator on $L^p(R)$, for all $1 \leq p \leq \infty$. The needed assumption turns out to be the same as the conditions considered in Sections 4 and 6. We also classify the WH-systems for which the frame operator extends to a bounded linear operator on the Wiener amalgam-space $W(L^\infty, \ell^1)$. Finally, in Section 8, we classify the WH-systems whose frame operators are equal. Again, this will be a natural consequence of the properties considered throughout the paper.

The first author would like to thank Nakhle Asmar, Nigel J. Kalton, and Stephen Montgomery-Smith for many helpful discussions during the preparation of this paper.

## 2. Basic Tools

A sequence $(f_i)_{i \in I}$ in a Hilbert space $H$ is a **frame** for $H$ if there are $A > 0$ and $B < \infty$ such that

$$(2.1) \qquad A\|f\|^2 \leq \sum_{k \in I} |<f, f_i>|^2 \leq B\|f\|^2, \quad \text{for all } f \in H.$$

We call $A$ the **lower frame bound** and $B$ the **upper frame bound** for the frame. If $A = B$, we say $(f_i)$ is a **tight frame** and if $A = B = 1$ it is called a **normalized tight frame**. A bounded unconditional basis for a Hilbert space $H$ is called a **Riesz basis** (See the end of this section for a detailed discussion of this topic). If $(f_i)$ and $(g_i)$ are frames for Hilbert spaces $H, K$ respectively, we say that $(f_i)$ is **equivalent to** $(g_i)$ is there is a bounded, linear invertible operator $L : H \to K$ satisfying $Lf_i = g_i$, for all $i \in I$. In particular, if $(f_i)$ is a frame for $H$ and $L$ is a bounded linear invertible operator from $H$ to $H$, then $(Lf_i)$ is also a frame for $H$ which is equivalent to $(f_i)$.

If $(e_i)_{i \in I}$ is an orthonormal basis for $\ell_2$ and $(f_i)_{i \in I}$ is a sequence of vectors in $H$, we will call the operator $T : \ell_2 \to H$ given by: $T(e_i) = f_i$ the **pre-frame operator** associated with $(f_i)_{i \in I}$. A direct calculation shows that

$$T^* f = \sum_{i \in I} <f, f_i> e_i, \quad \text{for all } f \in H.$$

It follows that the pre-frame operator is bounded if and only if there is a constant $B > 0$ satisfying:

$$\sum_{i \in I} |<f, f_i>|^2 \leq B\|f\|^2, \quad \text{for all } f \in H.$$

Moreover, $(f_i)_{i \in I}$ is a frame if and only if the pre-frame operator is a bounded, linear, onto map. In that case $S = TT^*$ is a invertible operator on $H$, called the **frame operator**, which has the form:

$$(2.2) \qquad Sf = \sum_{i \in I} <f, f_i> f_i, \text{ for all } f \in H.$$

Now, for all $f \in H$ we have,

$$(2.3) \qquad <Sf, f> = \sum_{i \in I} |<f, f_i>|^2.$$

So, the frame operator is a positive, self-adjoint invertible operator on $H$. It is a straightforward calculation that for all $f \in H$,

(2.4)
$$f = \sum_{i \in I} <S^{-1}f, f_i> f_i = \sum_{i \in I} <f, S^{-1}f_i> f_i = \sum_{i \in I} <f, S^{-1/2}f_i> S^{-1/2}f_i.$$

Since $S$ is an invertible operator on $H$, $(S^{-1}f_i)$ is a frame for $H$ which is equivalent to the frame $(f_i)$ and is called the **minimal dual frame** for $(f_i)$. In general, there may be other "dual frames" for $(f_i)$ in the sense that there may be sequences $(h_i)$ of elements of $H$ so that for every $f \in H$ we have

$$f = \sum_i <f, h_i> f_i, \text{ for all } f \in H.$$

It can be shown [13] that a bounded frame $(f_i)$ has only one dual frame if and only if $(f_i)$ is a Riesz basis for $H$. Also, [13] no two distinct dual frames can be equivalent.

To define Weyl-Heisenberg frames, let $a, b \in R$ and define the operators of **modulation** $E_b$, and **translation** $T_a$ for functions $f \in L^2(R)$ by:

$$E_b f(t) = e^{2\pi i b t} f(t),$$

and

$$T_a f(t) = f(t - a).$$

Given $g \in L^2(R)$, and $a, b > 0$, if the family $(E_{mb} T_{na} g)_{m,n \in Z}$ is a frame for $L^2(R)$, we call this a **Gabor frame** or a **Weyl-Heisenberg frame** or just a **WH-frame** for short. The function $g$ is referred to as the **window function**. The numbers $a, b$ are the **frame parameters**, with $a$ being the **shift parameter** and $b$ the **modulation parameter**. It is known that if $ab > 1$, then for any $g \in L^2(R)$, the family $(E_{mb}T_{na}g)$ is not a frame, and in fact, it is not even complete in $L^2(R)$ (see [4] p. 978). Also, if $ab = 1$, this family is a frame if and only if it is a Riesz basis for $L^2(R)$. Finally, if $ab < 1$, the family is **overcomplete** in the sense that one can remove at least one element from the frame and the remaining elements will still form a frame (but perhaps with different frame bounds).

Let $a, b > 0$. We will be interested in the frame operator for functions $g$ for which $(E_{mb}T_{na}g)_{m,n \in Z}$ has at least a finite upper frame bound. This class of functions is called the class of **preframe functions** and is denoted by PF. Thus we suppress the dependence on $a, b$. We have immediately from the definitions and our earlier considerations

**Proposition 2.1.** *The following are equivalent:*

*(1) $g \in PF$.*

*(2) The operator*

$$Sf = \sum_{n,m \in Z} <f, E_{mb}T_{na}g> E_{mb}T_{na}g,$$

*is a well defined bounded linear operator on $L^2(R)$.*

A **WH-system** is any $(g, a, b)$ with $g \in L^2(R)$, and a **preframe WH-system** is any $(g, a, b)$ with $g \in PF$.

If $ab$ is rational, we will call $(g, a, b)$ a **rational WH-system**. Wexler and Raz [25] have shown that for any two Weyl-Heisenberg systems $(E_{mb}T_{na}g)_{m,n \in Z}$ and $(E_{mb}T_{na}h)_{m,n \in Z}$ these systems are dual if and only if

(\*) $$<h, E_{k/b}T_{\ell/a}g> = ab\delta_{k0}\delta_{\ell 0}, \quad k, \ell \in Z,$$

where $\delta$ is the Kronecker delta. This gives rise to a natural dual frame, viz. that $h$ satisfying (\*) with minimal norm, called the **Wexler-Raz dual.** A direct calculation shows that the frame operator for a Weyl-Heisenberg frame commutes with both translation and modulation. In particular, the minimal dual frame of a Weyl-Heisenberg frame $(E_{mb}T_{na}g)_{m,n \in Z}$ with frame operator $S$ is itself a Weyl-Heisenberg frame of the form $(E_{mb}T_{na}S^{-1}g)_{m,n \in Z}$. The fact that the minimal dual of a Weyl-Heisenberg frame equals the Wexler-Raz dual is a rather subtle point. We refer the reader to [17], Propositions 3.2-3.3 and [6], Proposition 4.2.

For a general shift-invariant system $g_{nm} = (g_m(\cdot - na))_{n,m \in Z}$, we also have the notion of an **upper frame bound**. This is a constant $B < \infty$ such that

$$\sum_{n,m \in Z} |<f, g_{nm}>|^2 \leq B\|f\|^2, \quad f \in L^2(R).$$

There is a simple necessary and sufficient condition for the existence of finite upper frame bounds for shift-invariant systems [15], Sec. 1.2, [22], Sec. 1.4..

**Proposition 2.2.** *For a shift-invariant system $(g_m(\cdot - na))_{n,m \in Z}$, the following are equivalent:*

*(1) The system has a finite upper frame bound $B$.*

*(2) We have*

$$(2.5) \qquad \left( \frac{1}{a} \sum_m \hat{g}_m(\nu - k/a)\overline{\hat{g}_m(\nu - \ell/a)} \right)_{k,\ell \in Z} \leq B \cdot I.$$

*where $\hat{g}(\nu) = \int e^{-2\pi it\nu} g(t) dt$ denotes the Fourier transform of $g$.*

If we fix the row with index 0 in (2.5) above, we get a necessary condition for the existence of finite upper frame bounds.

**Proposition 2.3.** *If a shift-invariant system $(g_{mn})$ has a finite upper frame bound $B$ then*

$$(2.6) \qquad \sum_{k \in Z} |\frac{1}{a} \sum_{m \in Z} \hat{g}_m(\nu)\overline{\hat{g}_m(\nu - k/a)}|^2 \leq B \quad a.e.$$

*Proof.* The row with index 0 in (2.5) is just

$$\left( \frac{1}{a} \sum_{m \in Z} \hat{g}_m(\nu)\overline{\hat{g}_m(\nu - k/a)} \right)_{k \in Z}.$$

We will make extensive use of the WH-frame Identity due to Daubechies [4] (see also [14]). Careful examination of the proof of this result in [14], Theorem 4.1.5, shows that we do not need the assumption of $g \in$ PF used there. By checking functions which are bounded and supported just on an interval of length $1/b$, one can see that the proof goes through without change for any $g \in L^2(R)$ with the property that $\sum_n |g(x - na)|^2 \leq B$ a.e. It is this stronger result we will need and now state.

**WH-Frame Identity.** *If $\sum_n |g(t-na)|^2 \leq B$ a.e. and $f \in L^2(R)$ is bounded and compactly supported, then*

$$(2.7) \qquad \sum_{n \in Z} \sum_{m \in Z} |<f, E_{mb}T_{na}g>|^2 = F_1(f) + F_2(f)$$

where

$$(2.8) \quad F_1(f) = b^{-1} \int_R |f(t)|^2 \sum_n |g(t-na)|^2 \, dt,$$

$$(2.9) \quad F_2(f) = b^{-1} \sum_{k \neq 0} \int_R \overline{f(t)} f(t-k/b) \sum_n g(t-na) \overline{g(t-na-k/b)} \, dt =$$

$$b^{-1} \sum_{k \geq 1} 2 \operatorname{Re} \int_R \overline{f(t)} f(t-k/b) \sum_n g(t-na) \overline{g(t-na-k/b)} \, dt.$$

To simplify the notation a little we introduce the following auxiliary functions:

$$(2.10) \quad G_k(t) = \sum_{n \in Z} g(t-na) \overline{g(t-na-k/b)}, \quad \text{for all} \ \ k \in Z.$$

When it is necessary to keep track of the function $g$ in (2.10), we write $G_k \equiv G_{g,k}$. It follows easily from the WH-frame Identity that if $g \in \mathrm{PF}$, then there is a constant $B > 0$ such that,

$$(2.11) \quad G_0(t) = \sum_{n \in Z} |g(t-na)|^2 \leq B, \quad \text{a.e.}$$

Note that the $G_k$ are periodic functions on $R$ of period a.

Next we relate the condition (2.5) with some of the operators we will use throughout this paper.

**Proposition 2.4.** *Let $a, b \in R$ with $ab \leq 1$ and $g \in L^2(R)$ and assume that*

$$\sum_{k \in Z} |G_k(t)|^2 \leq B, \quad \text{a.e.}$$

*Then for all bounded, compactly supported functions $f \in L^2(R)$ the series*

$$Lf = b^{-1} \sum_k (T_{k/b} f) G_k,$$

*converges unconditionally in norm in $L^2(R)$. Moreover,*

$$<Lf, f> = \sum_{m,n \in Z} |<f, E_{mb} T_{na} g>|^2.$$

*Finally, if $g \in PF$, so that the series*

$$Sf = \sum_{m,n} <f, E_{mb}T_{na}g> E_{mb}T_{na}g,$$

*also converges unconditionally in $L^2(R)$, we have that $Lf = Sf$.*

*Proof.* First we will check that the series for $Lf$ is unconditionally convergent for all bounded, compactly supported functions $f \in L^2(R)$. Since this class of functions are finite sums of bounded functions supported on intervals of the form $I_n = [na, (n+1)a), n = 1, 2, \cdots$, it suffices to assume that $f$ is a bounded function supported on $I_n$ with uniform upper bound $D$. Now, since $a \leq 1/b$, we have that the functions $[(T_{k/b}f)G_k]_{k \in Z}$ are disjointly supported. Since the $G_k$ are periodic of period a, a simple calculation yields for all $M \subset Z$ with $|M| < \infty$,

$$\|\sum_{k \in M}(T_{k/b}f)G_k\|_{L^2(R)} = \int_0^a |f(t)|^2 \sum_{k \in M}|G_k(t)|^2 dt \leq D^2 \int_0^a \sum_{k \in M}|G_k(t)|^2 dt.$$

Since $\sum_k |G_k(t)|^2 \leq B$ a.e., it follows from the Monotone Convergence Theorem that the series defining $Lf$ (indeed, even every subseries) converges. So the series for $Lf$ converges unconditionally in $L^2(R)$.

For the moreover part, we check

$$<Lf, f> = <b^{-1}\sum_{k \in Z}(T_{k/b}f)G_k, f> = b^{-1}\sum_{k \in Z}<(T_{k/b}f)G_k, f>$$

$$= b^{-1}\sum_{k \in Z}\int_R \overline{f(t)}f(t-k/b)G_k(t)dt = \sum_{m,n \in Z}|<f, E_{mb}T_{na}g>|^2,$$

where the last equality follows from the WH-frame identity.

To see that $Sf = Lf$, we mimic the argument of [14],Theorem 4.2.1, to discover that for all bounded, compactly supported $f \in L^2(R)$, and for all $h \in L^2(R)$ we have

$$<b^{-1}\sum_k(T_{k/b}f) \cdot G_k, h> = \sum_{m,n}<f, E_{mb}T_{na}g><E_{mb}T_{na}g, h> = <Sf, h>.$$

We define the **Wiener amalgam space** $W(L^\infty, \ell^1)$ to be the set of all measurable functions $g$ on $R$ for which there is some $a > 0$ such that

$$\|g\|_{W,a} = \sum_{n \in Z}\|g \cdot \chi_{[an,a(n+1))}\|_\infty = \sum_{n \in Z}\|T_{na}g \cdot \chi_{[0,a)}\|_\infty < \infty.$$

It is easily checked that $W(L^\infty, \ell^1)$ is a Banach space with the norm $\|\cdot\|_{W,a}$. For the proof of the following properties of this space we refer to [14], Proposition 4.1.7.

**Lemma 2.5.** *For a function $g \in W(L^\infty, \ell^1)$ we have*
  (1) *If $\|g\|_{W,a}$ is finite for one $a$ then it is finite for all $a$.*
  (2) *If $m$ is a natural number and $0 < b \leq ma$, then $\|g\|_{W,a} \leq 2m\|g\|_{W,b}$.*

We end this section with a discussion concerning the types of convergence we will work with in this paper. If $x_n$ are elements of some Banach space $X$, a series $\sum_n x_n$ is said to be **unconditionally convergent** if for every increasing sequence of natural numbers $(k_n)$ we have that

$$\text{norm } \lim_n \sum_{j=1}^n x_{k_j}$$

exists.

The next result is well known and can be found for example in [20], Proposition 1.c.1.

**Proposition 2.6.** *For $x_n$ in a Banach space $X$, the following are equivalent:*
  (1) *$\sum_n x_n$ is unconditionally convergent.*
  (2) *$\sum_n x_{\sigma(n)}$ converges for every permutation $\sigma$.*
  (3) *$\sum_n \theta_n x_n$ converges for every choice of $\theta_n = 0, \pm 1$. (Or equivalently, for every choice of complex $|\theta_n| \leq 1$).*

*Moreover, in this case there is a constant $K$ so that for every choice of scalars $(a_n)$ we have*

$$\|\sum_n a_n x_n\| \leq K \sup_n |a_n| \|\sum_n x_n\|.$$

We will also work with weak convergence of series in $L^2(R)$. The celebrated Orlicz-Pettis theorem says that weak unconditional convergence and norm unconditional convergence are the same in every Banach space. This result can be fund in [7], Chapter IV, p. 24.

**Orlicz-Pettis Theorem.** *If $x_n$ are elements of a Banach space so that for every increasing sequence of natural numbers $(k_n)$ we have*

$$\text{weak } \lim_n \sum_{j=1}^n x_{k_j}$$

*exists, then the series $\sum_n x_n$ is unconditionally convergent.*

The problem in applications of this theorem is that we must know the vector which is the unconditional sum of the series. The notion of weakly unconditionally Cauchy, allows us to check unconditional convergence of a series just from the series itself without having to know what it is converging to.

**Definition 2.7.** *A series $\sum_n x_n$ is said to be **weakly unconditionally Cauchy** (wuC) if given any permutation $\sigma$ of the natural numbers we have that $(\sum_{k=1}^n x_{\sigma(k)})$ is a weakly Cauchy sequence.*

Recall the Banach space $c_0$:

$$c_0 = \{x = (a_n) : \|x\| =: \sup_n |a_n| < \infty \text{ and } \lim_{n \to \infty} a_n = 0\}$$

Another well-known result ([7], Section V, Theorem 6) gives criteria for checking when a series is wuC.

**Theorem 2.8.** *The following are equivalent for a series $\sum_n x_n$ in a Banach space:*

*(1) $\sum_n x_n$ is wuC.*
*(2) For every $x^* \in X^*$ we have*

$$\sum_n |x^*(x_n)| < \infty.$$

*(3) For every $(a_n) \in c_0$ we have that*

$$\sum_n a_n x_n$$

*converges.*

*(4) There is a constant $C$ so that for every finite subset $M$ of the natural numbers we have*

$$\|\sum_{n \in M} x_n\| \leq C.$$

Our last result of this section, combined with the above, shows that all these notions of unconditional convergence are equivalent for a Hilbert space. We can, and will, therefore use them interchangeably throughout the paper. First, we recall that a Banach space $Y$ is said to **embed** into a Banach space $X$ if there is a subspace $Z$ of $X$ which is isomorphic to $Y$. This theorem can be found in [7], P. 45, Theorem 8).

**Theorem 2.9.** *If $c_0$ does not embed into a Banach space $X$, then every series $\sum_n x_n$ which is wuC is also unconditionally convergent in $X$.*

## 3. The Zak Transform and Necessary Conditions for Finite Upper Frame Bounds

We will make extensive use of the Zak transform throughout the paper. So we introduce the basic properties here. For a good survey of this important subject see [14] Section 1.5, or [18].

**Definition 3.1.** *The **Zak transform** of a function $f \in L^2(R)$ is*

(3.1) $$(Z_\lambda f)(t, \nu) = \lambda^{1/2} \sum_{k \in Z} f(\lambda(t-k))e^{2\pi i k \nu}, \quad a.e. \ t, \nu \in R,$$

*where the right-hand side has to be interpreted in $L^2_{loc}(R^2)$-sense.*

In this section we will work with the case $\lambda = 1$, but in later sections we will also need $\lambda = a, 1/b$ etc. When $\lambda = 1$ we just write $Z$ for the (standard) Zak transform (omitting the subscript 1). So we will introduce the basic properties of the general Zak transform (see [13]) which we list below, where $f, h \in L^2(R)$:

(3.2) $$\lambda^{1/2} f(\lambda t) = \int_0^1 (Z_\lambda f)(t, \nu) d\nu, \quad a.e. \ t \in R.$$

(3.3) $$(Z_\lambda f)(t+1, \nu) = e^{2\pi i \nu}(Z_\lambda f)(t, \nu), \quad a.e. \ t, \nu \in R.$$

(3.4) $$(Z_\lambda f)(t, \nu+1) = (Z_\lambda f)(t, \nu), \quad a.e. \ t, \nu \in R.$$

(3.5) $$(Z_\lambda \hat{f})(t, \nu) = e^{2\pi i \nu t}(Z_{1/\lambda} f)(-\nu, t), \quad a.e. \ t, \nu \in R.$$

(3.6) $$\int_{-\infty}^\infty f(t)\overline{h(t)}dt = \int_0^1 \int_0^1 (Z_\lambda f)(t, \nu)\overline{(Z_\lambda h)(t, \nu)}dt d\nu.$$

Equalities (3.3) and (3.4), also referred to as the *quasi-periodicity relations*, imply that the Zak transform is completely determined by its values in the unit square $Q = [0, 1) \times [0, 1)$. Moreover, by (3.6), the Zak transform is a unitary map of $L^2(R)$ onto $L^2(Q)$.

**Proposition 3.2.** *Suppose that $(E_m T_n g)$ has a finite upper frame bound $B$. For each $k \in Z$ with $G_k$ defined by (2.10) we have*
  *(1) $|(Zg)(t, \nu)|^2 \leq B$, a.e. $t, \nu \in [0, 1)$.*

(2) $|(Zg)(t,\nu)|^2$ has for a.e. t the Fourier series expansion

$$|(Zg)(t,\nu)|^2 = \sum_{k=-\infty}^{\infty} G_k(t) e^{-2\pi i k \nu}, \quad a.e. \quad \nu.$$

(3) For a.e. t and all $k \in Z$,

$$G_k(t) = \int_0^1 |(Zg)(t,\nu)|^2 e^{2\pi i k \nu} d\nu.$$

*Proof.* (1) follows from the definition of the Zak transform. (2) is an application of Parseval's Identity and Carleson's theorem, and (3) follows immediately from (2).

For the next result, we need a simple lemma.

**Lemma 3.3.** *Let $F(t,\nu)$ be measurable in $(t,\nu) \in Q$, and assume that $|F(t,\nu)| \leq B$ a.e. $t,\nu$. Also assume that $F(t,\nu)$ is continuous in $\nu$ for a.e. t. Then,*

$$ess\ sup_t |F(t,\nu_0| \leq B, \quad for\ ALL\ \ \nu_0 \in [0,1).$$

*Proof.* Take $\nu_0 \in [0,1)$. We have for a.e. $t \in [0,1)$ that

(3.7) $$|F(t,\nu)| \leq B, \quad a.e. \ \ \nu \in [0,1).$$

Take a $t \in [0,1)$ such that $F(t,\cdot)$ is continuous and (3.7) holds for this t. Then we can find a sequence $(\nu_n)$ in [0,1) such that $\nu_n \to \nu_0$, while $|F(t,\nu_n)| \leq B$. By the continuity of $F(t,\nu)$ in $\nu$, it follows that $|F(t,\nu_0| \leq B$. Hence, ess $sup_t |F(t,\nu_0)| \leq B$, as required.

Now we give the basic construction for our examples. We will refine this considerably in Sections 5 and 6.

**Proposition 3.4.** *Let $g \in L^2(R)$ and assume that $(E_m T_n g)$ has a finite upper frame bound. For $M \subset Z$ with $|M| < \infty$ and all $f \in L^2(R)$ let $S_M : L^2(R) \to L^2(R)$ be given by*

$$S_M f = \sum_{k \in M} f(\cdot - k) G_k = \sum_{k \in M} (T_k f) \cdot G_k.$$

*Then*

$$\|S_M\| = ess\ sup_{t,\nu} |\sum_{k \in M} G_k(t) e^{-2\pi i k \nu}|.$$

Moreover, for all $\nu_0 \in [0, 1)$ we have

$$\text{ess sup}_t |\sum_{k \in M} G_k(t) e^{-2\pi i k \nu_0}| \leq \|S_M\|.$$

*Proof.* For any $f \in L^2(R)$, taking the Zak transform, we have

$$Z \left[ \sum_{k \in M} f(\cdot - k) G_k \right] (t, \nu) = (Zf)(t, \nu) \cdot \sum_{k \in M} G_k(t) e^{-2\pi i k \nu}, \quad \text{a.e. } t, \nu \in R.$$

Since $Z$ is unitary, the operator norm of $S_M$ is the same as the operator norm of the multiplication operator:

$$Zf \in L^2(Q) \to (Zf)(t, \nu) \cdot \sum_{k \in M} G_k(t) e^{-2\pi i k \nu}.$$

But these operator norms are precisely,

$$\|S_M\| = \text{ess sup}_{t,\nu} |\sum_{k \in M} G_k(t) e^{-2\pi i k \nu}|.$$

For the moreover part of the proposition, we just take

$$F(t, \nu) = \sum_{k \in M} G_k(t) e^{-2\pi i k \nu}$$

in Lemma 3.3.

We end this section by using the Zak transform to examine some special consequences of a WH-system having a finite upper frame bound when a=b=1.

**Proposition 3.5.** *If $(g, 1, 1)$ has a finite upper frame bound $B$ then for all $n = 1, 2, \cdots$ and for all $\nu_0 \in [0, 1)$*

$$G_0(t) + 2 \sum_{k=1}^{n} \frac{n-k}{n} Re \left[ G_k(t) e^{-2\pi i k \nu_0} \right] \leq B, \quad \text{a.e. } t.$$

*Proof.* If $(g, 1, 1)$ has upper frame bound B then by Proposition 2.3 we have

$$\sum_{n=-\infty}^{\infty} |g(t-n)|^2 \leq B,$$

and by Proposition 3.2 we have

$$|(Zg)(t,\nu)|^2 \leq B, \quad \text{a.e.} \quad \nu \in [0,1).$$

Since $G_{-k}(t) = \overline{G_k(t)}$, we have for any $\nu_0$ that

$$G_0(t) + 2\sum_{k=1}^{n} \frac{n-k}{n} Re\left[G_k(t)e^{-2\pi ik\nu_0}\right]$$

$$= \sum_{k=-n}^{n}\left(1-\frac{|k|}{n}\right)G_k(t)e^{-2\pi ik\nu_0} = \int_0^1 |(Zg)(t,\nu)|^2 \sum_{k=-n}^{n}\left(1-\frac{|k|}{n}\right)e^{2\pi ik(\nu-\nu_0)}d\nu$$

$$= \int_0^1 |(Zg)(t,\nu)|^2 \frac{1}{n}\left(\frac{\sin\ \pi n(\nu-\nu_0)}{\sin\ \pi(\nu-\nu_0)}\right)^2 d\nu \leq B\int_0^1 \frac{1}{n}\left(\frac{\sin\ \pi n(\nu-\nu_0)}{\sin\ \pi(\nu-\nu_0)}\right)^2 d\nu = B.$$

This is Cesaro summation of a Fourier series, so we get for all $\nu_0$ that

$$G_0(t) + 2\sum_{k=1}^{n} \frac{n-k}{n} Re\left[G_k(t)e^{-2\pi ik\nu_0}\right] \leq B, \quad \text{a.e.} \quad t.$$

The condition given in Proposition 3.5 is actually quite delicate. We will now show that we cannot replace the numbers $(n-k)/n$ by ones in Proposition 3.5.

**Example 3.6.** *There is a WH-system $(g,1,1)$ with a finite upper frame bound for which*

$$\sum_{k=1}^{\infty} Re\left[G_k(t)\right] = \infty$$

*This, of course, is the condition in Proposition 3.5 with (n-k)/n replaced by ones and $\nu_0 = 1$.*

*Proof.* In [9], page 161-162, Edwards constructs a continuous, even, real periodic f (which we can clearly take of period 1) such that the partial sums

(3.8) $$\sum_{k=-n}^{n} a_k e^{-2\pi ik\nu} \quad \text{with} \quad a_k = \int_0^1 f(\nu)e^{2\pi ik\nu}d\nu,$$

are unbounded at $\nu = 0$. Now take $M \in R$ such that $M + f(\nu) > 0$ everywhere, and let $g$ be such that

$$|(Zg)(t,\nu)|^2 = M + f(\nu), \quad (t,\nu) \in Q.$$

Since $f$ is even and real, we have that the $G_k(t) = G_k$ (independent of t) are real and even in $k \in Z$ as well. But these $G_k$ equal the $a_k$ of (3.8) above for $k \neq 0$. Letting $\nu_0 = 0$, we have the desired example.

In the special case of Proposition 3.5 when $g$ is positive and real valued, Proposition 3.5 yields a necessary condition for $(E_m T_n g)$ to have a finite upper frame bound. In general, this condition is called the CC-condition and will be examined in detail in section 4.

**Corollary 3.7.** *Let $g \in L^2(R)$ be positive and real valued. Then the following are equivalent:*
  *(1) The WH-system $(g, 1, 1)$ has a finite upper frame bound.*
  *(2) There is a constant $K > 0$ such that*

$$\sum_{k \in Z} |G_k(t)| \leq K, \quad a.e. \ t \in R.$$

*Proof.* (1) $\Rightarrow$ (2): Letting $\nu_0 = 0$ in Proposition 3.5, we see that our asumptions imply there is a constant $B > 0$ so that

$$\sup_n \sum_{k=1}^{n} \frac{n-k}{n} G_k(t) \leq B, \quad a.e. \ t.$$

Now (2) follows easily from here.

(2) $\Rightarrow$ (1): This is a result of Casazza and Christensen [2]

The next example shows that the positivity assumption in Corollary 3.7 is necessary.

**Example 3.8.** *There is a real valued function $g$ such that $(E_m T_n g)$ has a finite upper frame bound, while $\sum_k |G_k|$ is not essentially bounded.*

*Proof.* We define

$$F(t, \nu) = \begin{cases} 1 & : 0 \leq t < 1, \ 0 \leq \nu < 1/4 \\ \frac{1}{2} & : 0 \leq t < 1, \ \frac{1}{4} \leq \nu < \frac{3}{4} \\ 1 & : 0 \leq t < 1, \ \frac{3}{4} \leq \nu < 1 \end{cases}$$

Now extend $F$ quasi-periodically by equations (3.3) and (3.4) to all of $R \times R$. Choose the unique function $g \in L^2(R)$ for which $F = Zg$. Using equality (3.2) we can find

this $g$ explicitly as

$$g(t+k) = \begin{cases} \dfrac{3}{4} & : k = 0, \ 0 \le t < 1 \\ \dfrac{\sin \frac{1}{2}\pi k}{2\pi k} & : k \ne 0, \ 0 \le t < 1 \end{cases}$$

A direct calculation (see e.g. [3], the proof of Proposition 2.4) shows that

$$\sum_k |G_k| = \infty, \quad \text{everywhere.}$$

## 4. Sufficient Conditions for Finite Upper Frame Bounds

It is well known [14] that a function $g$ in the Wiener space $W(L^\infty, \ell^1)$ is in PF. This was the standard condition used until recently for guaranteeing a function in $L^2(R)$ yielded a WH-system with a finite upper frame bound. Recently, Casazza and Christensen [2], gave a significantly weaker condition which guarantees a function is in PF.

**Theorem 4.1.** *Let $a, b \in R$ and $g \in L^2(R)$. If there is a constant $B > 0$ so that*

(CC) $$\sum_{k \in Z} |G_k(t)| \le B, \quad a.e.$$

*then $g \in PF$.*

After introduction in [2] of the condition in Theorem 4.1, this condition has become known in small circles as the CC-condition. To avoid one and the same condition having different names we shall adhere to this convention. In [3], there is an example of a function in PF which fails to satisfy the CC-condition above. It is, however, easy to see that this condition is much weaker than requiring that $g$ be in the Wiener space.

**Example 4.2.** *There are functions satisfying the CC-condition which are not in the Wiener space.*

*Proof.* We let $a = b = 1$ and let

$$g = \sum_{n=0}^{\infty} \chi_{[n+1-2^n, n+1-2^{n+1})}.$$

Since this $g$ is positive and real, we can calculate

$$\sum_{k \in Z} |\sum_{n \in Z} g(t-n)\overline{g(t-n-k)}| = \left(\sum_{n \in Z} g(t-n)\right)^2 = 1.$$

Clearly, $\sum_n \|g \cdot \chi_{[m,m+1)}\|_\infty = \infty$.

In this section we will give a detailed analysis of the CC-condition and its relation to the frame bounds. To begin, we will examine the connection between this condition and Zak matrices. We refer the reader to [15] Section 1.5 for more details. Let $ab = p/q$ with $p,q \in Z$ and gcd(p,q)=1. We make a particular choice for the Zak transform,

$$(4.1) \qquad (Zf)(t,\nu) = b^{-1/2} \sum_{\ell=-\infty}^{\infty} f(\frac{t-\ell}{b})e^{2\pi i \ell \nu}, \quad \text{a.e. } t,\nu \in R.$$

where $f \in L^2(R)$. We define for $f,h \in L^2(R)$, and a.e. $t,\nu \in R$, the $(p \times q)$-matrix,

$$(4.2) \qquad \Phi^f(t,\nu) = p^{-1/2}\left((Zf)(t-\ell\frac{p}{q}, \nu + \frac{k}{p})\right)_{k=0,1,\cdots,p-1;\ell=0,1,\cdots,q-1}$$

and the $(p \times p)$-matrix (where * denotes conjugate transpose)

$$(4.3) \qquad A^{fh}(t,\nu) = \Phi^f(t,\nu)[\Phi^h(t,\nu)]^*.$$

When $(E_{mb}T_{na}g)$ has a finite upper frame bound, then for $f \in L^2(R)$ and a.e. $t,\nu$ we have

$$(4.4) \qquad \Phi^{Sf}(t,\nu) = A^{gg}(t,\nu)\Phi^f(t,\nu).$$

Hence the action of the frame operator amounts to multiplication in the Zak-domain by the matrix $A^{gg}(t,\nu)$. In particular, $S$ has the frame bounds $A,B$ if and only if

$$(4.5) \qquad A \cdot I_{p \times p} \leq A^{gg}(t,\nu) \leq B \cdot I_{p \times p} \quad \text{a.e. } t,\nu.$$

We can now calculate for $f \in L^2(R)$, $n = 1, 2, \cdots$ and a.e. $t,\nu$

$$(4.6) \qquad \Phi^{S^n f}(t,\nu) = (A^{gg}(t,\nu))^n \Phi^f(t,\nu).$$

Also, denoting by ${}^\circ\gamma = S^{-1}g$, the minimal dual window of the WH-frame generated by $g$, we have for $f \in L^2(R)$ and a.e. $t,\nu$,

$$(4.7) \qquad A^{{}^\circ\gamma {}^\circ\gamma}(t,\nu) = [A^{gg}(t,\nu)]^{-1}$$

(4.8) $$\Phi^{S^{-1}f}(t,\nu) = A^{\circ\gamma\circ\gamma}(t,\nu)\Phi^f(t,\nu).$$

Formulas (4.6) - (4.8) are special cases of the formula

(4.9) $$\Phi^{\phi(S)f}(t,\nu) = \phi(A^{gg}(t,\nu))\Phi^f(t,\nu) \quad \text{a.e. } t,\nu$$

for $f \in L^2(R)$, where $\phi$ is a function analytic in a neighborhood of $[A, B]$. For some special cases like $\phi(x) = x^\alpha$, we can obtain this result directly by using

$$x^\alpha = \left(\frac{A+B}{2}\right)^\alpha \left(1 - (1 - \frac{2}{A+B}x)\right)^\alpha$$

$$= \left(\frac{A+B}{2}\right)^\alpha \sum_{n=0}^\infty \binom{\alpha}{n}(-1)^n \left(1 - \frac{2}{A+B}x\right)^n.$$

Note that $I - \frac{2}{A+B}S$ has operator norm $\leq \frac{B-A}{B+A} < 1$. For general $\phi$, we may have to use a Dunford representation

$$\phi(S) = \frac{1}{2\pi i}\int_C (\lambda I - S)^{-1}\phi(\lambda)d\lambda$$

(where C is a contour containing $[A, B]$ in its interior) together with (4.6) to make sense out of this calculation.

The connection between the CC-condition and the Zak matrices $A^{gg}(t,\nu)$ can now be obtained. From [15], Theorem 1.2.5 and Subsection 1.3.2, we have $g \in$ PF with upper frame bound $B$ if and only if for a.e. t,

$$\frac{1}{b}M_g(t)M_g^*(t) = \left(\frac{1}{b}\sum_n g(t - (\ell/b) + na)\overline{g(t - (j/b) + na)}\right)_{\ell,j \in Z} \leq B \cdot I.$$

Here, $M_g(t) = (g(t + na - \ell/b))_{\ell,n \in Z}$. Therefore, for a.e. t we have

(4.10) $$\left(\sum_n g(t + (m/b) - na)\overline{g(t + (m/b) - na - \ell/b)}\right)_{\ell \in Z} =$$

$$\left(\overline{M_g(t+m/b)M_g^*(t+m/b)e_0}\right)_{\ell \in Z} \in \ell^2,$$

with $\ell^2$-norm $\leq bB$ for all m and a.e. t. Thus, for all $m, k = 0, 1, \cdots, p-1$ and a.e. t

(4.11) $$S_{mk}^{gg}(t,\nu) =$$

$$\sum_\ell \left(\sum_n g(t + (m/b) - na)\overline{g(t + (m/b) - na - (\ell/b))}\right) \cdot e^{-2\pi i\ell(\nu + k/p)},$$

is well-defined as a 1-periodic $L^2[0,1)$-function in $\nu$. Since $M_g(t)M_g^*(t)$ and $M_g^*(t)M_g(t)$ have the same operator norm, we also have that for a.e. t,

(4.12) $$\sum_k |g(t-k/b)|^2 \leq bB.$$

We now have

**Proposition 4.3.** *Let $ab = p/q$ with $\gcd(p,q) = 1$, and let $g \in PF$. The following are equivalent:*

*(1) $g$ satisfies the CC-condition.*

*(2) $S^{gg}(t,\nu)$ has for a.e. t an absolutely convergent Fourier series in $\nu$ with a.e. t-independent bound on the sums of the absolute values of the Fourier coefficients.*

*Proof.* We observe that

$$\sum_n g(t+(m/b)-na)\overline{g(t+(m/b)-na-\ell/b)} = G_{m-\ell}(t+m/b).$$

Thus,

$$\sum_\ell |\sum_n g(t+(m/b)-na)\overline{g(t+(m/b)-na-\ell/b)}e^{-2\pi i\ell k/p}| =$$

$$\sum_\ell |G_{m-\ell}(t+m/b)|, \quad m,k = 0,1,\cdots,p-1.$$

This easily implies the equivalence of (1) and (2).

We can now relate $S^{gg}$ to $A^{gg}$.

**Proposition 4.4.** *If $g \in PF$, then*

$$S_{mk}^{gg}(t,\nu) = be^{-2\pi imk/p} \sum_{r=0}^{p-1} A_{rk}^{gg}(bt,\nu)e^{2\pi imr/p}, \quad a.e. \ t,\nu.$$

*and*

$$A_{j,k}^{gg}(bt,\nu) = \frac{1}{bp} \sum_{m=0}^{p-1} S_{m,k}^{gg}(t,\nu)e^{2\pi im(k-j)/p}, \quad a.e. \ t,\nu.$$

*Proof.* Assume that $g \in$ PF. We want to show that for a.e. $t,\nu$

(4.13) $$S_{mk}^{gg}(t,\nu) = be^{-2\pi imk/p} \sum_{r=0}^{p-1} A_{rk}^{gg}(bt,\nu)e^{2\pi imr/p},$$

and

$$(4.14) \quad A^{gg}_{jk}(t,\nu) = \frac{1}{bp} \sum_{m=0}^{p-1} S^{gg}_{mk}(t,\nu) e^{2\pi i(k-j)/p},$$

for $m, k, j = 0, 1, \cdots, p-1$.

Note that when $f \in L^2(R)$ we have for all t outside a null set $N$ that $\sum_\ell |f(\frac{t-\ell}{b})|^2$ is finite a.e. For such a t we have that

$$(4.15) \quad \sum_{\ell=-\infty}^{\infty} f(\frac{t-\ell}{b}) e^{2\pi i \ell \nu}, \quad \text{a.e.} \quad \nu \in R,$$

defines a 1-periodic $L^2[0,1)$-function of $\nu$, and we shall choose one representative of $Zf$ in (4.2) such that $(Zf)(t,\nu)$ agrees with (4.1) a.e. $\nu \in R$ when t is outside N. For $g \in$ PF, the situation is even slightly better as (4.12) shows.

Now we let t be such that $bt - j\frac{p}{q} \notin N$ for $j = 0, 1, \cdots, q-1$, and we shall compute the Fourier coefficients with respect to $\nu$ of the 1-periodic function

$$(4.16) \quad A^{gg}_{rk}(t,\nu) = \frac{1}{p} \sum_{j=0}^{q-1} (Zg)(bt - j\frac{p}{q}, \nu + \frac{r}{p}) \overline{(Zg)(bt - j\frac{p}{q}, \nu + \frac{k}{p})}.$$

Let $\ell \in Z$. There holds (since $\frac{p}{qb} = a$)

$$(4.17) \quad \int_0^1 e^{2\pi i\ell\nu} (Zg)(bt - j\frac{p}{q}, \nu + \frac{r}{p}) \overline{(Zg)(bt - j\frac{p}{q}, \nu + \frac{k}{p})} d\nu =$$

$$\frac{1}{b} \int_0^1 e^{2\pi i\ell\nu} (\sum_{s=-\infty}^{\infty} g(t - ja - \frac{s}{b}) e^{2\pi is(\nu + \frac{r}{p})}) \overline{(\sum_{j=-\infty}^{\infty} g(t - ja - \frac{s}{b}) e^{2\pi is(\nu + \frac{k}{p})})} d\nu$$

$$= \frac{1}{b} \int_0^1 (\sum_{s=-\infty}^{\infty} g(t - ja - \frac{s}{b}) e^{2\pi is\frac{r}{p}}) \overline{(\sum_{s=-\infty}^{\infty} g(t - ja - \frac{s}{b}) e^{2\pi is\frac{r}{p}} e^{2\pi is\nu})} d\nu$$

$$= \frac{1}{b} \int_0^1 (\sum_{s=-\infty}^{\infty} g(t - ja - \frac{s-\ell}{b}) e^{2\pi i(s+\ell)\nu} e^{2\pi is\nu}) \overline{(\sum_{s=-\infty}^{\infty} g(t - ja - \frac{s}{b}) e^{2\pi is\frac{-k}{p}} e^{2\pi is\nu})} d\nu$$

$$= \frac{1}{b} \sum_{s=-\infty}^{\infty} g(t - ja - \frac{s-\ell}{b}) \overline{g(t - ja - \frac{s}{b})} e^{2\pi i(s-\ell)\frac{r}{p} - 2\pi is\frac{k}{p}},$$

where for the last equality we have used Parseval's formula. Thus,

$$(4.18) \quad \int_0^1 e^{2\pi i\ell\nu} \sum_{r=0}^{p-1} A^{gg}_{rk}(bt,\nu) e^{2\pi im\frac{r}{p}} d\nu =$$

$$\frac{1}{bp}\sum_{j=0}^{q-1}\sum_{r=0}^{p-1}\sum_{s=-\infty}^{\infty} g(t-ja-\frac{s-\ell}{b})\overline{g(t-ja-\frac{s}{b})}e^{-2\pi is\frac{k}{p}}e^{2\pi i(s-\ell+m)\frac{r}{p}}.$$

Carrying out the summation over r in (4.18), we see that only s of the form $\ell - m + vp$, $v \in Z$ survive. And thus (4.18) is equal to

$$\frac{1}{b}\sum_{j=0}^{q-1}\sum_{v=-\infty}^{\infty} g(t-ja-\frac{-m+vp}{b})\overline{g(t-ja-\frac{\ell-m-vp}{b})}e^{-2\pi i(\ell-m)\frac{k}{p}} =$$

$$\frac{1}{b}\sum_{j=0}^{q-1}\sum_{v=-\infty}^{\infty} g(t-(ja+\frac{vp}{b})+\frac{m}{b})\overline{g(t-(ja+\frac{vp}{b})+\frac{m}{b}-\frac{\ell}{b})}e^{-2\pi i(\ell-m)\frac{k}{p}}.$$

Now
$$ja + \frac{vp}{b} = (j + \frac{vp}{ab})a = (j+vq)a,$$

and $j+vq$, ranges through all of Z when $j = 0, 1, \cdots, q-1$, $v \in Z$. Hence, (4.18) is equal to

$$\frac{1}{b}\sum_{j=-\infty}^{\infty} g(t+\frac{m}{b}-ja)\overline{g(t+\frac{m}{b}-ja-\frac{\ell}{b})}e^{-2\pi i(\ell-m)\frac{k}{p}}.$$

Therefore,
$$\int_0^1 be^{-2\pi im\frac{k}{p}}\sum_{r=0}^{p-1} A^{gg}_{rk}(bt,\nu)e^{2\pi im\frac{r}{p}}e^{2\pi i\ell\nu}d\nu =$$

$$e^{-2\pi i\ell\frac{k}{p}}\sum_{j=-\infty}^{\infty} g(t+\frac{m}{b}-ja)\overline{g(t+\frac{m}{b}-ja-\frac{\ell}{b})}.$$

We now choose t also such that (4.10) holds. For such a t we then have that

$$be^{-2\pi im\frac{k}{p}}\sum_{r=0}^{p-1} A^{gg}_{rk}(bt,\nu)e^{2\pi im\frac{r}{p}}$$

$$\sim \sum_{\ell=-\infty}^{\infty}\left(\sum_{j=-\infty}^{\infty} g(t+\frac{m}{b}-ja)\overline{g(t+\frac{m}{b}-ja-\frac{\ell}{b})}\right)e^{-2\pi i\ell\frac{k}{p}}e^{-2\pi i\ell\nu} = S^{gg}_{mk}(t,\nu).$$

This establishes (4.13). The identity (4.14) follows directly from (4.13), and this completes the proof of Proposition 4.4.

Note that both $S^{gg}(bt,\nu)$ and $A^{gg}(t,\nu)$ are periodic in $t,\nu$ with periods a,1 respectively. By Proposition 4.4, it follows that for a particular t, $A^{gg}(t,\nu)$ has an absolutely convergent Fourier series in $\nu$ if and only if $S^{gg}(t,\nu)$ has such a series.

**Remark.** *In Proposition 4.4 and the results following Definition 4.7, we could have worked with the Zak transform $Z_a$ instead of $Z_{1/b}$ and we would have obtained the same results. This is based on [15], Subsec. 1.5.7 and the fact that for all $k, \ell = 0, 1, \ldots, p - 1$ and a.e. $t, \nu$*

$$\frac{1}{p}\sum_{j=0}^{q-1}(Z_a g)(t - k\frac{q}{p}, \nu - \frac{j}{q})(Z_a g)^*(t - \ell\frac{q}{p}, \nu - \frac{j}{q}) = \frac{1}{b}\sum_{r=-\infty}^{\infty} G_{\ell-k+rp}(at - \frac{k}{b})e^{-2\pi i r q \nu},$$

*when $g \in \mathbf{PF}$. These developments are discussed in detail at the end of this section.*

We can now easily pass the CC-condition through the frame operator even for general shift-invariant systems.

**Theorem 4.5.** *Let $(g_m(\cdot - na))_{m,n \in Z}$ be a shift-invariant system with a finite upper frame bound and frame operator $S$. If the system satisfies the CC-condition, then $((Sg)(\cdot - na))_{n,m \in Z}$ satisfies the CC-condition.*

*Proof.* Let $B$ be the upper frame bound of our system. We have for $f \in L^2(R)$ and a.e. $\nu \in R$ that
$$(\hat{S}f)(\nu) = \sum_{\ell} \hat{S}_{o\ell}(\nu)\hat{f}(\nu - \ell/a),$$

where
$$\hat{S}_{0\ell}(\nu) = \frac{1}{a}\sum_{m} \hat{g}_m(\nu)\overline{\hat{g}_m(\nu - \ell/a)}.$$

This can also be written in terms of the matrix

$$\hat{S}(\nu) = \left(\frac{1}{a}\sum_{m}\hat{g}_m(\nu - k/a)\overline{\hat{g}_m(\nu - \ell/a)}\right)_{k,\ell \in Z}$$

as
$$\left((\hat{S}f)(\nu - k/a)\right)_{k \in Z} = \hat{S}(\nu)(\hat{f}(\nu - \ell/a))_{\ell \in Z}.$$

¿From this it is evident that we have for polynomials $\phi$ that

$$\left((\phi(\hat{S})f)(\nu - k/a)\right)_{k \in Z} = \phi(\hat{S}(\nu))(\hat{f}(\nu - \ell/a))_{\ell \in Z}$$

and one can even allow functions $\phi$ analytic in an open set contained in the closed segment $[A, B]$. This is so since the frame operator $S$ has the frame bounds $A, B$ if and only if $AI \leq \hat{S}(\nu) \leq BI$ a.e. $\nu \in R$.

Now assume that our system satisfies the CC-condition. That is, we assume that we have an $M > 0$ such that for a.e. $\nu \in R$

(4.19) $$\sum_\ell |\frac{1}{a} \sum_m \hat{g}_m(\nu)\hat{g}_m(\nu - \ell/a)| le M.$$

Since $S$ commutes with all shifts $T_{na}$, $n \in Z$, it follows that the system $((Sg_m)(\cdot - na))_{n,m \in Z}$ has frame operator $S^3$. Then by the above functional calculus we must show that there is a $K > 0$ such that for a.e. $\nu \in R$

$$\sum_\ell |\hat{S}^3_{0\ell}(\nu)| \leq K.$$

Usine (4.19) with $\nu - k/a$ instead of $\nu$, we see that there is an $M > 0$ such that for a.e. $\nu \in R$

$$\sum_\ell |\hat{S}_{k\ell}(\nu)| \leq M.$$

Therefore, for all $i \in Z$ and a.e. $\nu \in R$

$$\sum_j |\hat{S}^2_{ij}(\nu)| = \sum_j |\sum_k \hat{S}_{ik}(\nu)\hat{S}_{kj}(\nu)| \leq$$

$$\sum_k |\hat{S}_{ik}(\nu)| \sum_j |\hat{S}_{kj}(\nu)| \leq M^2,$$

and, for all $i \in Z$ and a.e. $\nu \in R$

$$\sum_\ell |\hat{S}^3_{i\ell}(\nu)| = \sum_\ell |\sum_k \hat{S}^2_{ik}(\nu)\hat{S}_{k\ell}(\nu)| \leq$$

$$\sum_k |\hat{S}^2_{ik}(\nu)| \sum_\ell |\hat{S}(\nu)| \leq M^3.$$

This is what we wanted (even somewhat more).

It is an open question whether the Wexler-Raz dual $S^{-1}g$ satisfies the CC-condition if $g$ satisfies it. Later in this section we will examine a "uniform" version of the CC-condition and show that this condition is indeed inherited by $S^{-1}g$.

We will now show that the CC-condition for shift invariant systems yields a finite upper frame bound for these systems. In this setting, if

$$(g_{mn}) = (g_m(\cdot - na))_{n,m \in Z}$$

is a shift invariant system, the CC-condition looks like:

(4.20) $$\sum_{k \in Z} |\frac{1}{a} \sum_m \hat{g}_m(\nu)\overline{\hat{g}_m(\nu - k/a)}| \leq B, \quad \text{a.e. } \nu.$$

First, we recall the Schur Test.

**Schur Test.** *Assume that $A = (A_{k,\ell})_{k,\ell \in Z}$ satisfies*

$$\sum_{k \in Z} |A_{k,\ell}| \leq B, \quad \ell \in Z \quad \text{and} \quad \sum_{\ell \in Z} |A_{k,\ell}| \leq B, \quad k \in Z.$$

*Then $A$ defines a bounded linear operator of $\ell^2$ with $\|A\|_{\ell^2 \to \ell^2} \leq B$.*

Our next result now follows easily.

**Proposition 4.6.** *If a shift invariant system $(g_m(\cdot - na))_{m,n \in Z}$ satisfies the CC-condition, then it has a finite upper frame bound.*

*Proof.* If we assume the CC-condition, since countable unions of null sets are null sets, we have for a.e. $\nu$ that

$$\sum_k |\frac{1}{a} \sum_m \hat{g}_m(\nu - k/a) \overline{\hat{g}_m(\nu - \ell/a)}| \leq B, \quad \forall \ell \in Z,$$

and

$$\sum_\ell |\frac{1}{a} \sum_m \hat{g}_m(\nu - k/a) \overline{\hat{g}_m(\nu - \ell/a)}| \leq B, \quad \forall k \in Z.$$

By the Schur Test, the hypotheses of Proposition 2.2 are satisfied, and so the system has a finite upper frame bound.

Now we will consider a uniform version of the the CC-condition which is strong enough to pass to the dual.

**Definition 4.7.** *Assume that $(g, a, b)$ has a finite upper frame bound. We say that $(g, a, b)$ satisfies a **uniform CC-condition** when for every $\epsilon > 0$ there is a $K > 0$ so that for a.e. $t$ we have*

(UCC) $$\sum_{|k| \geq K} |G_k(t)| = \sum_{|k| \geq K} |\sum_n g(t - na) \overline{g(t - na - k/b)}| < \epsilon.$$

It is easily checked that $g$ satisfies the UCC-condition if $g \in W(L^\infty, \ell^1)$. Also, using the finite upper frame bound assumption, it is easily checked that the UCC-condition implies the CC-condition. The converse implication fails as the next example shows.

**Example 4.8.** *There is a WH-system $(g, 1, 1)$ satisfying the CC-condition but failing the UCC-condition.*

*Proof.* In the case $a = b = 1$ and $g \in L^2(R)$, the CC-condition means there is a $B > 0$ such that for a.e. t,

$$\sum_{k \in Z} |\sum_{n \in Z} g(t - n) \overline{g(t - n - k)}| \leq B,$$

while the UCC-condition means for all $\epsilon > 0$ there is a $K > 0$ such that for a.e. $t$,

$$\sum_{k \in Z, \ |k| \geq K} |\sum_{n \in Z} g(t-n)\overline{g(t-n-k)}| < \epsilon.$$

Also, we note that for a $g$ satisfying the CC (or UCC)-condition, the quantities

$$G_k(t) = \sum_{n \in Z} g(t-n)\overline{g(t-n-k)}, \quad k \in Z$$

are the Fourier coefficients of $|(Zg)(t,\nu)|^2$ (see Section 3). Now, choose real numbers $\alpha_\ell > 0$, for $\ell = 0, 1, 2 \cdots$ with

$$M = \sum_{\ell=0}^{\infty} \alpha_\ell < \infty,$$

and define $F(t,\nu)$ for $t, \nu \in [0,1)$ by setting

$$F(t,\nu) = M + \sum_{\ell=m}^{\infty} \alpha_{\ell-m} \cos 2\pi\ell\nu, \quad t \in [1 - \frac{1}{m}, 1 - \frac{1}{m+1}), \quad \nu \in [0,1).$$

Then let $g$ be the unique element of $L^2(R)$ with

$$(Zg)(t,\nu) = \sqrt{F(t,\nu)}, \quad t, \nu \in [0,1).$$

Now, a direct calculation shows that for this $g$ we have

$$\sum_n g(t-n)\overline{g(t-n-k)} = \begin{cases} 0 & : \ |k| = 0, 1, \cdots, m-1, \\ \frac{1}{2}\alpha_{|k|-m} & : \ |k| = m, m+1, \cdots \end{cases}$$

when $t \in [1 - \frac{1}{m}, 1 - \frac{1}{m+1})$. Hence, $g$ satisfies the CC-condition but not the UCC-condition.

We can also formulate the UCC-condition in terms of the behavior of the Fourier series for $A^{gg}(t, \cdot)$.

**Proposition 4.9.** *For $ab = \frac{p}{q}$, with gcd $(p,q) = 1$, and $g \in PF$, the following are equivalent*

*(1) $g$ satisfies the uniform CC-condition.*

(2) there is a null set $N$ of $t$'s such that the Fourier series of $A^{gg}(t, \cdot)$ converges uniformly outside $N$.

*Proof.* Assume $g \in \text{PF}$ satisfies the UCC-condition. For $j = 1, 2, \cdots$ we let $K_j > 0$ and the null set $N_j \subset R$ be such that

$$\sum_{|\ell| \geq K_j - p + 1} |G_\ell(t)| \leq \frac{1}{j}, \quad t \notin N_j.$$

We let $N = \cup_{j=1}^\infty N_j$ and finally let $W = \cup_{m=0}^{p-1}(N - \frac{m}{b})$. Outside the null set $W$ we have

$$\sum_{|\ell| \geq K_j} |\sum_n g(t + \frac{m}{b} - na)\overline{g(t + \frac{m}{b} - na - \frac{\ell}{b})} e^{-2\pi i \ell \frac{k}{p}}| =$$

$$\sum_{|\ell| \geq K_j} |G_{m-\ell}(t + \frac{m}{b})| \leq \sum_{|\ell| \geq K_j - p + 1} |G_\ell(t + \frac{m}{b})| \leq \frac{1}{j},$$

for $m, k = 0, 1, \cdots, p - 1$. Thus we see that the Fourier series for $S_{mk}^{gg}(t, \cdot)$ is uniformly convergent for $t \notin W$, $m, k = 0, 1, \cdots, p - 1$. Evidently, from Proposition 4.4, the same holds for the $A_{jk}^{gg}(t, \cdot)$.

The proof of the converse is similar.

To establish that the uniform CC-condition is really a usable tool, we next show that a condition of Tolimieri and Orr [15], Subsec. 1.4.3, is strong enough to imply it.

**Definition 4.10.** *Let $(g, a, b)$ be a WH-system with finite upper frame bound. For $x, y \in R$, write for reasons of conciseness*

$$g_{x,y}(t) = e^{2\pi i y t} g(t - x) = (E_y T_x g)(t), \quad t \in R.$$

*We say that $(g, a, b)$ satisfies the condition of Tolimieri and Orr when*

(A) $$\sum_{k, \ell} | < g, g_{k/b, \ell/a} > | < \infty.$$

It is known [15], Subsec. 1.4.3, that condition A guarantees unconditional convergence of the Janssen representation for the frame operator:

$$Sf = \frac{1}{ab} \sum_{k, \ell} < g, g_{k/b, \ell/a} > f_{k/b, \ell/a}, \quad f \in L^2(R).$$

We can also formulate Condition A in terms of the Zak matrices (see [19], Theorem 1.3).

**Proposition 4.11.** *Let $g \in PF$. Then the following are equivalent*

*(1) $g$ satisfies Condition A.*

*(2) The matrix $A^{gg}(t, \nu)$ (which is periodic in both variables) has an absolutely convergent Fourier series in both variables.*

Now we relate Condition A to the CC-condition.

**Proposition 4.12.** *If $(g, a, b)$ is a WH-system with a finite upper frame bound, and $g$ satisfies condition A, then $g$ satisfies the uniform CC-condition.*

*Proof.* Since $(g, a, b)$ has a finite upper frame bound, by equality (2.11) there is a $B > 0$ so that

$$(4.21) \qquad \sum_n |g(t - na)|^2 \leq B, \quad \text{a.e. } t.$$

Hence, the functions

$$G_k(t) = \sum_n g(t - na)\overline{g(t - na - k/b)}, \quad k \in Z$$

are well-defined a.e. with finite $L^\infty$-norm. Note that $G_k(t)$ is periodic in t with period a. By (4.5) we can compute the Fourier coefficients of $G_k(t)$ by integrating term by term according to

$$\frac{1}{a}\int_0^a G_k(t) e^{-2\pi i \ell t/a} dt = \frac{1}{a}\int_0^a \sum_n g(t - na)\overline{g(t - na - k/b)} e^{-2\pi i \ell (t - na)/a} dt$$

$$= \frac{1}{a}\int_{-\infty}^\infty g(t)\overline{g(t - k/b)} e^{-2\pi i \ell t/a} dt = \frac{1}{a} <g, g_{k/b, \ell/a}> .$$

Thus,

$$(4.22) \qquad G_k(t) \sim \frac{1}{a}\sum_\ell <g, g_{k/b, \ell/a}> e^{-2\pi i \ell t/a}.$$

Now, if we assume condition A then the right-hand side of (4.22) is a continuous function and agrees a.e. with $G_k(t)$. Thus,

$$|G_k(t)| \leq \frac{1}{a}\sum_\ell |<g, g_{k/b, \ell/a}>|, \quad \text{a.e. } t.$$

and so

$$\sum_k |G_k(t)| \leq \frac{1}{a}\sum_{k,\ell} |<g, g_{k/b, \ell/a}>| < \infty, \quad \text{a.e. } t.$$

Now it is clear that the uniform CC-condition is satisfied.

The next example shows that the converse of Proposition 4.12 fails in general.

**Example 4.13.** *There is a function $g \in PF$ satisfying the uniform CC-condition but not satisfy condition A.*

*Proof.* We let $a = b = 1$ and take

$$g(t) = \chi_{[0,1/2)}(t) + \frac{1}{2}\chi_{[1/2,1)}(t).$$

Then

$$(Zg)(t,\nu) = g(t) \text{ and } A^{gg}(t,\nu) = g^2(t), \quad 0 \leq t,\nu < 1,$$

are discontinuous. Whence, $A^{gg}$ cannot possess an absolutely convergent Fourier series in two variables. By Proposition 4.11, it follows that $g$ does not satisfy Condition A. However, for any $t$, $A^{gg}(t,\cdot)$ has a Fourier series in $\nu$ consisting of the constant term only which is therefore uniformly absolutely convergent with respect to t. So $g$ satisfies the uniform CC-condition by Proposition 4.9.

Finally, we will show that the uniform CC-condition passes to the Wexler-Raz dual $S^{-1}g$. We first need a uniform version of a theorem of Wiener. For completeness, we will sketch a proof.

**Uniform Wiener $\frac{1}{f}$-Theorem.** *Let $V$ be a set of 1-periodic functions with uniformly absolutely convergent Fourier series. That is, for every $\epsilon > 0$, there is a $K \in N$ such that*

$$(4.23) \qquad \sum_{|k| \geq K} |a_k(f)| \leq \epsilon, \quad \text{for all } f \in V,$$

*where*

$$a_k(f) = \int_0^1 e^{-2\pi i k \nu} f(\nu) d\nu.$$

*Also, assume there are real numbers $0 < C, D$ such that*

$$(4.24) \qquad C \leq |f(\nu)| \leq D, \quad \text{for all } f \in V, \; \nu \in R.$$

*Then for every $\epsilon > 0$, there is a $K > 0$ so that*

$$(4.25) \qquad \sum_{|k| \geq K} |a_k(1/f)| < \epsilon, \quad \text{for all } f \in V.$$

*Proof.* We follow the elementary proof of Wiener's 1/f-Theorem due to Newman [21]. Accordingly, we can assume that

$$|f(\nu| \geq 1, \quad \text{for all } f \in V, \; \nu \in R.$$

Also, by (4.23) and (4.24) and some elementary considerations, we have that there is a real number $E > 0$ such that

$$\|f\| = \sum_{k=-\infty}^{\infty} |a_k(f)| \leq E, \quad \text{for all } f \in V.$$

Again by (4.23) and (4.24), we can find a $K \in N$, an $F > 0$ and for any $f \in V$, a trigonometric polynomial $P_f$ of degree $\leq K$ such that

$$\|P_f - f\| \leq \frac{1}{3} \quad \text{and} \quad \max|P_f| \leq F.$$

Now by Bernstein's Theorem ([1], Theorem 11.1.2) we have that

$$\max|P_f'| \leq 2\pi K \max|P_f| \leq 2\pi K F, \quad \text{for all } f \in V.$$

We thus see that all the bounds that play a role in Newman's proof hold uniformly for $f \in V$, and hence the result follows.

**Theorem 4.14.** *If ab is rational, $(g, a, b)$ a Weyl-Heisenberg frame with frame operator $S$ and if $g$ satisfies the uniform CC-condition, then $S^{-1}g$ satisfies the uniform CC-condition.*

*Proof.* Let $°\gamma = S^{-1}g$. To prove the theorem, we check that outside a null set of t's, we have uniformly absolutely convergent Fourier series for all matrix elements $S_{mk}^{°\gamma°\gamma}(t, \cdot)$ as a function of the second variable. Then the result will follow from Proposition 4.9.

We let $N_1$ be a null set of t's outside which we have

$$A \leq \frac{1}{b} \sum_{n} |g(t - na)|^2 \leq B,$$

and

$$\sum_{|k| \geq K_m} |\sum_{n} g(t - na)\overline{g(t - na - k/b)}| \leq \frac{1}{m},$$

for some $0 < A, B < \infty$, $k_1 < k_2, \cdots$. By Proposition 4.9, all matrix elements $S_{mk}^{gg}(t, \cdot)$, whence all matrix elements $A_{rk}^{gg}(bt, \cdot)$, have uniformly absolutely convergent Fourier series in the second variable, for t outside $N_1$. Since $(g, a, b)$ is a WH-frame, we know by (4.5) that a.e. $t, \nu$

(4.26) $$AI_{p \times p} \leq A^{gg}(bt, \nu) \leq BI_{p \times p}.$$

So we can find a null set $N_2$ of t's such that for t outside $N_2$ we have that (4.26) holds a.e. in $\nu \in R$. Letting $N = N_1 \cup N_2$, we have for t outside $N$ that all matrix elements of $A^{gg}(bt, \nu)$ have absolutely and uniformly convergent Fourier series in $\nu$ and are continuous in $\nu$. Therefore, (4.26) holds for all $\nu$. It follows that for all t outside $N$, we have that $\det[A^{gg}(bt, \nu)]$ has an absolutely uniformly convergent Fourier series in $\nu$ and is uniformly bounded away from 0 and $\infty$.

We now consider for t outside $N$ (see (4.7)),

$$A^{\circ\gamma\circ\gamma}(bt, \nu) = \frac{\mathrm{adj}(A^{gg}(bt, \nu))}{\det[A^{gg}(bt, \nu)]},$$

where we have used Cramer's Rule for the computation of the inverse of the matrix. By the uniform $\frac{1}{f}$-Theorem, and some elementary properties of absolutely convergent Fourier series, we conclude that all matrix elements of $A^{\circ\gamma\circ\gamma}(bt, \nu)$ have absolutely and uniformly convergent Fourier series in $\nu$, whenever $t \notin N$. This implies, by Proposition 4.4, that all matrix elements of $S^{\circ\gamma\circ\gamma}(t, \nu)$ have absolutely and uniformly convergent Fourier series in $\nu$ whenever $b \notin N$. Now by Proposition 4.9, the uniform CC-condition is satisfied by $(^\circ\gamma, a, b)$.

As noted by Newman [21] at the end of his proof, we can allow more general functions $\phi$ analytic in a neighborhood of the range of $f$. The uniform $\frac{1}{f}$-Theorem can then be generalized accordingly. Hence, one is led to expect validity of the uniform CC-condition for systems like $(S^\alpha g, a, b)$ for certain $\alpha$ whenever $(g, a, b)$ is a frame satisfying the uniform CC-condition. We have not carried out this program in detail, but it would be interesting if a careful study was made along these lines. The proof cannot, however, be based on Cramer's rule as we used in the proof of Theorem 4.14 above, but should rely now on Dunford representations

$$(A^{gg}(bt, \nu))^\alpha = \frac{1}{2\pi i} \int_C (zI - A^{gg}(bt, \nu))^{-1} z^\alpha dz,$$

together with the uniform $\frac{1}{f}$-Theorem.

When we want to replace in Theorem 4.14 the uniform CC-condition by the CC-condition, we are faced with the fact that there does not seem to be available a Wiener $\frac{1}{f}$-Theorem that yields a bound on $\|1/f\|$ in terms of the bounds $M, C, D$ bounding $f$ by: $\|f\| \leq M$, and $0 < C \leq |f(\nu)| \leq D < \infty$. This also seems to be an interesting direction of study in itself.

We end this section by presenting an alternative approach to Proposition 4.4 and the results following Definition 4.7 which use the Zak transform $Z_a$ in the place of

$Z_{1/b}$. Let for $f, h \in L^2(R)$ and a.e. $t, \nu$ (see [15], Subsec. 1.5.7)

(4.28) $$\Psi^f(t,\nu) = \left(p^{-1/2}(Z_a f)(t - k\frac{q}{p}, \nu - \frac{j}{q})\right)_{k=0,\ldots,p-1; j=0,\ldots,q-1},$$

(4.29) $$B^{fh}(t,\nu) = \Psi^f(t,\nu)\overline{\Psi^h(t,\nu)}.$$

Then for $g \in \mathbf{PF}$ and $f \in L^2(R)$ there holds a.e. $t, \nu$

(4.30) $$\Psi^{Sf}(t,\nu) = B^{gg}(t,\nu)\Psi^f(t,\nu),$$

and $B^{gg}$ is 1-periodic in t and $q^{-1}$-periodic in $\nu$. Now, Tolimieri-Orr's Condition A is satisfied if and only if $B^{gg}$ has an absolutely convergent Fourier series. Next, in place of Proposition 4.4, we have

**Proposition 4.15.** *We have for $k, \ell = 0, 1, \ldots, p-1$ and a.e. $t, \nu$*

(4.31) $$(B^{gg}(t,\nu))_{k,\ell} = \frac{1}{b}\sum_{r=-\infty}^{\infty} G_{\ell-k+rp}(at - \frac{k}{b})e^{-2\pi i r q \nu}.$$

*Proof.* We compute for $k, \ell = 0, 1, \ldots p-1$, $n \in Z$ and a.e. $t, \nu$

(4.32) $$\int_0^1 (B^{gg}(t,\nu))_{k,\ell} e^{2\pi i n \nu} d\nu =$$

$$\frac{1}{p}\int_0^1 \sum_{j=0}^{q-1}(Z_a g)(t - k\frac{q}{p}, \nu - \frac{j}{q})\overline{(Z_a g)(t - \ell\frac{q}{p}, \nu - \frac{j}{q})} e^{2\pi i n \nu} d\nu$$

$$= \frac{a}{p}\sum_{j=0}^{q-1}\int_0^1 \left[e^{2\pi i n \nu}\left(\sum_{s=-\infty}^{\infty} g(a(t - k\frac{q}{p} - s))e^{2\pi i s(\nu - j/q)}\right)\right.$$

$$\left.\left(\sum_{s=-\infty}^{\infty} \overline{g(a(t - \ell\frac{q}{p} - s))e^{2\pi i s(\nu - j/q)}}\right)\right] d\nu.$$

We move the $e^{2\pi i n \nu}$ to the second series over $s$ and then change the summation variable $s$ into $s + n$ to obtain

(4.33) $$\int_0^1 (B^{gg}(t,\nu))_{k,\ell} e^{2\pi i n \nu} d\nu$$

$$= \frac{a}{p}\sum_{j=0}^{q-1}\int_0^1 \left[\left(\sum_{s=-\infty}^{\infty} g(a(t - k\frac{q}{p} - s))e^{-2\pi i s j/q} \cdot e^{2\pi i s \nu}\right)\right.$$

$$\left(\sum_{s=-\infty}^{\infty}\overline{g(a(t-\ell\frac{q}{p}-s-n))}e^{-2\pi i(s+n)j/q}\cdot e^{2\pi is\nu}\right)\Bigg]d\nu$$

$$=\frac{a}{p}\sum_{j=0}^{q-1}\sum_{s=-\infty}^{\infty}g(a(t-k\frac{q}{p}-s))e^{-2\pi isj/q}\overline{g(a(t-\ell\frac{q}{p}-s-n))}e^{-2\pi i(s+n)j/q},$$

where in the last step we have used Parseval's formula. Now, since $\sum_{j=0}^{q-1}e^{2\pi inj/q}=q$ when $n$ is a multiple of $q$ and 0 otherwise, we obtain

(4.34) $$\int_0^1 (B^{gg}(t,\nu))_{k,\ell}e^{2\pi in\nu}d\nu =$$

$$\begin{cases} 0 & , \ n \neq 0 \pmod{q} \\ \frac{aq}{p}\sum_{s=-\infty}^{\infty}g(a(t-k\frac{q}{p}-s))\overline{g(a(t-\ell\frac{q}{p}-s-rq))}, & n=rq. \end{cases}$$

Now, since $\frac{aq}{p}=\frac{1}{b}$, (4.31) follows from (4.33), (4.34) and

$$\sum_{s=-\infty}^{\infty}g(a(t-k\frac{q}{p}-s))\overline{g(a(t-\ell\frac{q}{p}-s-rq))}=$$

$$\sum_{s=-\infty}^{\infty}g(at-\frac{k}{b}-sa)\overline{g(at-\frac{\ell}{b}-sa-\frac{rp}{b})}=$$

$$\sum_{s=-\infty}^{\infty}g(at-\frac{k}{b}-sa)\overline{g(at-\frac{k}{b}-sa-\frac{rp+\ell-k}{b})}=G_{rp+\ell-k}(at-\frac{k}{b}).$$

It follows that, for a particular t, all $(B^{gg}(t,\nu))_{k,\ell}$, $k,\ell=0,1,\ldots,p-1$, have an absolutely convergent Fourier series in $\nu$ if and only if

$$\sum_r |G_{rp+\ell-k}(at-\frac{k}{b})| < \infty, \quad k,\ell=0,1,\ldots,p-1.$$

¿From this point onwards, the developments after Proposition 4.4 can be essentially copied from what appears there.

Although in the preceding developments, we really had two different approaches depending upon our choice of the Zak transform, there are several situations later in the paper wither the choice of Zak transform is more or less dictated by the situation. For example, in the proof of Theorem 6.8, (2) $\Rightarrow$ (1), the choice for $Z_a$ is more or less dictated by the fact that the $G_k$ are a-periodic. In the proof of Proposition 6.11, the choice for $Z_{1/b}$ is also more or less dictated since choosing $Z_a$ instead would lead to a formula in which the dependence on k would also appear in the $Z_a f$.

# 5. Symmetric, Weak and Norm Convergence of the Walnut Representation

In this section we will make a detailed analysis of the convergence properties of the Walnut representation of the WH-frame operator generated by a function $g \in$ PF. Recall from Proposition 2.4 that for $g \in$ PF we have that the two series

$$(5.1) \qquad Lf = b^{-1} \sum_{k \in Z} (T_{k/b}f)G_k = Sf = \sum_{m,n} <f, E_{mb}T_{na}g> E_{mb}T_{na}g$$

converge unconditionally in norm on a dense subset of $L^2(R)$.

We will need some notation for checking convergence of series of arbitrary functions in the left-hand side of (5.1).

**Definition 5.1.** *Let $g \in L^2(R)$ satisfies $G_0(t) \leq B$ a.e. For any $f \in L^2(R)$, and any $K, L \in Z$, we let*

$$S_{K,L}f(t) = \sum_{k=-L}^{K} f(t - k/b) G_k(t),$$

*where as usual $G_k(t) = \sum_n g(t-na)\overline{g(t-na-k/b)}$. We also let $S_K = S_{K,K}$. If $M \subset Z$ with $|M| < \infty$, define*

$$S_M f(t) = \sum_{k \in M} f(t - k/b) G_k(t).$$

*Also, if $\mathcal{F}$ is the Fourier transform, let $T_{K,L} = \mathcal{F} S_{K,L}$. That is, for all $h \in L^2(R)$ and $\nu \in R$,*

$$(T_{K,L}h)(\nu) = \mathcal{F}(S_{K,L}h)(\nu).$$

If $\lim_{K \to \infty} S_K f$ exists, we say that **the Walnut series for $f$ converges symmetrically** - and this can be in either the norm or the weak topology - and we say **the Walnut series for f converges** when $\lim_{K,L \to \infty} S_{K,L}f$ exists. Our first result shows that weak and norm symmetric convergence of the Walnut series are equivalent.

**Theorem 5.2.** *For $a, b \in R$ and $g \in L^2(R)$ with $|G_0(t)| \leq B$ a.e., the following are equivalent:*
*(1) The Walnut series converges in norm symmetrically for every $f \in L^2(R)$.*
*(2) The Walnut series converges weakly symmetrically for every $f \in L^2(R)$.*

(3) We have $sup_K \|S_K\| = B < \infty$.

Moreover, in this case the WH-system $(g, a, b)$ has a finite upper frame bound and the Walnut series converges symmetrically to $Sf$, for all $f \in L^2(R)$.

*Proof.* (1) $\Rightarrow$ (2): this is obvious.

(2) $\Rightarrow$ (3): This is obvious when we interpret it correctly. Our assumption (2) implies that for every $f \in L^2(R)$, $(S_K f)_{K \in Z}$ converges weakly in $L^2(R)$, and hence is a bounded sequence. Since the $S_K$ are clearly bounded operators, the Uniform Boundedness Principle yields (3).

(3) $\Rightarrow$ (1): By (5.1) we know that $(S_K f)$ converges on the dense set of compactly supported $f \in L^2(R)$. By assumption (3), the operators $(S_K)$ are uniformly bounded. So $(S_K f)$ converges in norm for all $f \in L^2(R)$.

For the moreover part of the theorem, we note that (3) and Proposition 2.4 yield that there is a constant $B > 0$ so that for all bounded, compactly supported $f \in L^2(R)$ we have

$$\sum_{m,n \in Z} |<f, E_{mb}T_{na}g>|^2 = <Lf, f> \leq \|L\|\|f\|^2 \leq B\|f\|^2.$$

So our WH-system has upper frame bound $B$ for a dense set of elements of $L^2(R)$, whence for all of $L^2(R)$.

**Proposition 5.3.** *If the WH-system $(g, 1, 1)$ has a finite upper frame bound, then the following are equivalent:*

*(1) The Walnut series converges symmetrically for every $f \in L^2(R)$.*

*(2) There is a $B > 0$ such that*

$$sup_K \ ess \ sup_{t,\nu} |\sum_{k=-K}^{K} G_k(t) e^{-2\pi i k \nu}| \leq B.$$

*Proof.* By Proposition 3.4 (with $M = \{-K, -(K-1), \cdots, K-1, K\}$) we know that

$$\|S_K\| = ess \ sup_{t,\nu} |\sum_{k=-K}^{K} G_k(t) e^{-2\pi i k \nu}|.$$

Combining this with Theorem 5.2 yields the result.

Proposition 5.3, as well as Propositions 5.6 and 6.2 below have generalizations to the case of rational WH-systems. We will discuss these generalizations at the end of Section 6.

Our next example shows that there are functions $g \in$ PF for which the Walnut representation does not converge symmetrically for some function $f \in L^2(R)$.

**Example 5.4.** *There is a WH-system $(g, a, b)$, with a finite upper frame bound, for which the Walnut series does not converge symmetrically in $L^2(R)$, for some $f \in L^2(R)$.*

*Proof.* We let $a = b = 1$. We are going to construct a $g \in L^2(R)$ (actually, the Zak transform of such a $g$) which takes on the constant value $g_n$ on the intervals $[n, n+1)$. Note that for such a $g$, for each $k \in Z$ we have,

$$\rho_k = \sum_{n \in Z} g(t-n)\overline{g(t-n-k)} = \sum_{n \in Z} g_n \overline{g(n-k)},$$

are well-defined absolutely convergent series. Moreover, by Proposition 3.2, these numbers $\rho_k$ are the Fourier coefficients of $|(Zg)(t, \cdot)|^2$, when $0 \leq t < 1$. To construct our $g$, choose a continuous, real, 1-periodic function $\rho$, see [9], p. 161, so that

(5.2) $$\left( \sum_{k=-K}^{K} c_k e^{-2\pi i k \nu} \right)_{k=1,2,\ldots}$$

is an unbounded sequence at a particular point $\nu = \nu_0$, where the $c_k$ are the Fourier coefficients,

$$c_k = \int_0^1 e^{2\pi i k \nu} \rho(\nu) \, d\nu, \quad k \in Z.$$

Choose a real number $M$ so that $M + \rho(\nu) > 0$, for all $\nu \in [0, 1)$. Now define $g \in L^2(R)$ by specifying its Zak transform on $[0, 1) \times [0, 1)$ by

$$(Zg)(t, \nu) = (M + \rho(\nu))^{1/2}, \quad t, \nu \in [0, 1).$$

¿From our earlier discussion, we have for this $g$ that

$$\rho_k = M \delta_{k_0} + c_k, \quad k \in Z,$$

where $\delta_k$ is the Kronecker delta. Next, we will compute the norm of the operator $S_K$ for this $g$. Choose $f \in L^2(R)$ and for any $k = 1, 2, \cdots$ we consider

$$S_k f = \sum_{k=-K}^{K} \rho_k f(\cdot - k).$$

Proposition 3.4 yields,

$$\|S_K\| = \text{ess sup} \left| \left( M + \sum_{k=-K}^{K} c_k e^{2\pi i k \cdot} \right) \right|,$$

and so by (5.2),
$$\sup_K \|S_K\| = \infty.$$
By Theorem 5.2, we have that there is a $f \in L^2(R)$ so that $(S_K f)_{K \in Z}$ diverges.

Now we will examine the general convergence properties of the Walnut representation. We first show that weak and norm convergence of the Walnut series are equivalent properties.

**Theorem 5.5.** *For $a, b \in R$ and $g \in PF$, the following are equivalent:*
*(1) The Walnut series converges in norm for every $f \in L^2(R)$.*
*(2) The Walnut series converges weakly for every $f \in L^2(R)$.*
*(3) We have that $\sup_{K,L} \|S_{K,L}\| < \infty$.*
*Moreover, in this case the Walnut series converges to $Sf$, for all $f \in L^2(R)$.*

*Proof.* This is done exactly like the proof of Theorem 5.2 using Proposition 3.4. Again, the moreover part follows from Proposition 2.4.

Again we have a special version of Theorem 5.5 for a=b=1. We will not give a formal proof for this result since it follows exactly the format of the proof of Proposition 5.3.

**Proposition 5.6.** *If the WH-system $(g, 1, 1)$ has a finite upper frame bound, then the following are equivalent:*
*(1) The Walnut series converges weakly (or equivalently, strongly) for every $f \in L^2(R)$.*
*(2) There is a $B > 0$ such that for all $\nu_0 \in [0, 1)$,*
$$\sup_{K,L} \text{ess sup}_t \Big| \sum_{k=-L}^{K} G_k(t) e^{-2\pi i k \nu_0} \Big| \leq B,$$

At the end of Section 6, we generalize Proposition 5.6 to rational WH-systems.

A natural question is whether there is a connection between symmetric convergence of the Walnut series and norm convergence of the series. Our next example shows that these two types of convergence are really different.

**Example 5.7.** *For $a = b = 1$, there is a function $g \in PF$ for which the Walnut series converges symmetrically for every $f \in L^2(R)$, but there are functions $h \in L^2(R)$ for which the Walnut series does not converge.*

*Proof.* We proceed as in Example 5.3 but this time take a $\rho$ of the form:
$$\rho = \chi_{[0, 1/2)} + \frac{1}{2} \chi_{[1/2, 1)}.$$

For this $\rho$, the partial sums
$$\sum_{k=-K}^{K} \rho_k e^{-2\pi i k \nu}$$
are uniformly bounded and convergent a.e. to $\rho(\nu)$. For $f \in L^2(R)$, let $h = \mathcal{F}f$. Then we have
$$(T_{K,K}h)(\nu) = h(\nu) \sum_{k=-K}^{K} \rho_k e^{-2\pi i k \nu}.$$
Therefore, $T_{K,K}h \to h \cdot \rho$. It follows that the Walnut representation converges symmetrically for every $f \in L^2(R)$.

Arguing as in Proposition 2.4 of [3], we see that the $\rho_k$ are of the order of $1/k$, for each $k \in Z$. But, now
$$\sum_{k=1}^{K} \frac{1}{k} e^{-2\pi i k \nu} \sim \log K \quad \text{at} \quad \nu = 0.$$

By Proposition 5.6, there is an $h \in L^2(R)$ so that $(T_{K,L}h)$ is not bounded (and hence not convergent) in $L^2(R)$.

## 6. Unconditional Convergence of the Walnut Representation

We now address the question of when the Walnut series for the frame operator converges unconditionally. With the notation as in Definition 5.1, we first state the corresponding result of Theorem 5.5 for unconditional convergence.

**Theorem 6.1.** *Let $a, b \in R$ and $g \in PF$. The following are equivalent:*
  *(1) the Walnut series converges weakly unconditionally for every $f \in L^2(R)$.*
  *(2) The Walnut series converges unconditionally in norm for every $f \in L^2(R)$.*
  *(3) We have: $\sup_{M \subset Z, |M| < \infty} \|S_M\| < \infty$.*

*Proof.* The equivalence of (1) and (2) follows directly from the Orlicz-Pettis Theorem (see Section 2). The equivalence of (2) and (3) is proved in the same way as Theorems 5.2 and 5.5. That is, (3) is the statement that the bounded operators $S_M$ are uniformly bounded. But by Proposition 2.4, we know that these operators converge unconditionally for all $f$ in a dense subset of $L^2(R)$. So, they are uniformly bounded on $L^2(R)$ if and only if they converge pointwise (i.e. are pointwise bounded).

We again have the special case for a=b=1.

**Proposition 6.2.** *If the WH-system $(g, 1, 1)$ has a finite upper frame bound, then the following are equivalent:*

(1) *The Walnut series converges unconditionally for every $f \in L^2(R)$.*

(2) *There is a $B > 0$ such that for all $\nu_0 \in [0, 1)$,*

$$sup_M \, ess \, sup_t | \sum_{k \in M} G_k(t) e^{-2\pi i k \nu_0} | \leq B.$$

Proposition 6.2 will be generalized to rational WH-systems at the end of this section.

Before we examine unconditional convergence for the Walnut series, we need to note that this really is a different form of convergence than the norm convergence studied in Section 5.

**Example 6.3.** *There is a WH-system $(g, 1, 1)$ with $g \in PF$ for which the Walnut representation of the frame operator converges in norm for every $f \in L^2(R)$, but there is an $h \in L^2(R)$ so that the Walnut series for this $h$ does not converge unconditionally.*

*Proof.* In [27], there is a careful analysis of the series

$$(6.1) \qquad \sum_{k=1}^{\infty} k^{-\beta} e^{ik^\alpha} e^{-2\pi i k \nu},$$

with $0 < \alpha < 1$ and $\beta > 0$. In particular, it is shown that when $\beta > 1 - \frac{1}{2}\alpha$, then the series in (6.1) converges uniformly to a continuous sum $\psi_{\alpha,\beta}(-2\pi\nu)$. If we take say $\alpha = 1/2$ and $\beta = 7/8$, the condition is satisfied. Now, let the real continuous, 1-periodic function $\rho$ be given by,

$$\rho(\nu) = M + \sum_{k \neq 0} |k|^{-\beta} e^{i|k|^\alpha sgn(k)} e^{-2\pi i k \nu},$$

with $M > 0$ such that $\rho(\nu) > 0$ everywhere, Now we proceed as in Examples 5.3 and 5.5. If we choose the unique $g \in L^2(R)$ with

$$(Zg)(t, \nu) = (\rho(\nu))^{1/2}, \quad t, \nu \in [0, 1),$$

then

$$\left( \sum_{k=-L}^{K} \rho_k e^{-2\pi i k \nu} \right)_{L, K \in Z}$$

is uniformly bounded in the $L^\infty$-sense, and so (as in examples 5.3 and 5.5) we have that the Walnut series
$$\sum_{k \in Z} f(\cdot - k) G_k$$
converges in norm for every $f \in L^2(R)$. However, since $\sum_{k \in Z} |\rho_k|$ diverges, the Walnut series cannot converge unconditionally for every $f \in L^2(R)$ by Proposition 6.2.

Next we show that the CC-condition is strong enough to imply unconditional convergence of the Walnut representation.

**Theorem 6.4.** *If $G_k$ is a family of functions satisfying:*
*(1) $\sum_k |G_k| \leq B$.*
*(2) $\sum_k |T_{-k/b} G_k| \leq B$ a.e.*
*then for every $f \in L^2(R)$ the Walnut series*
$$\sum_k (T_{k/b} f) G_k$$
*converges unconditionally in $L^2$-norm.*

*Proof.* For any $h \in L^2(R)$ and any $m \geq n > 0$ we have:

$$(6.2) \quad \sum_{|k|=n}^{m} |<h, (T_{k/b}f)G_k>| = \sum_{|k|=n}^{m} |\int \overline{h(t)}(T_{k/b}f)(t)G_k(t)\, dt| \leq$$

$$\sum_{|k|=n}^{m} \int |h(t)||T_{k/b}f(t)||G_k(t)|\, dt = \sum_{|k|=n}^{m} \int |h(t)|\sqrt{|G_k(t)|}|T_{k/b}f(t)|\sqrt{|G_k(t)|}\, dt$$

$$\leq \sum_{|k|=n}^{m} (\int |h(t)|^2|G_k(t)|\, dt)^{1/2} (\int |T_{k/b}f(t)|^2|G_k(t)|\, dt)^{1/2} \leq$$

$$(\sum_{|k|=n}^{m} \int |h(t)|^2|G_k(t)|\, dt)^{1/2} (\sum_{|k|=n}^{m} \int |T_{k/b}f(t)|^2|G_k(t)|\, dt)^{1/2} =$$

$$(\int |h(t)|^2 \sum_{|k|=n}^{m} |G_k(t)|\, dt)^{1/2} (\int |f(t)|^2 \sum_{|k|=n}^{m} |T_{-k/b}G_k(t)|\, dt)^{1/2}.$$

By our hypotheses,
$$|h|^2 \sum_k |G_k| \in L^1(R) \quad \text{and} \quad |f|^2 \sum_k |T_{-k/b}G_k| \in L^1(R).$$

By the Lebesgue Dominated Convergence theorem we conclude that the following series converge in $L^1(R)$:

$$\sum_k |h(t)|^2 |G_k(t)| \text{ and } \sum_k |f(t)|^2 |T_{-k/b}G_k(t)|.$$

It follows that the right side of (6.2) goes to zero as $n \to \infty$. We conclude that the series

$$\sum_k (T_{k/b}f(t))G_k(t)$$

is weakly unconditionally Cauchy in $L^2(R)$. By Theorem 2.9, it is unconditionally convergent in norm.

For WH-systems we have that $T_{k/b}G_k = \overline{G_{-k}}$. So an immediate consequence of Theorem 6.4 is,

**Corollary 6.5.** *If $(g, a, b)$ satisfies the CC-condition, then for all $f \in L^2(R)$ the Walnut series $\sum_k (T_{k/b}f)G_k$ converges unconditionally.*

We have a converse for Corollary 6.4, even for more general systems, in the case where $ab$ is rational. But, before we present this, we give an example to show that things can become awkward if $ab$ is not rational.

**Example 6.6.** *There is irrational $ab$ and a function $g \in PF$, such that there is a $B > 0$ satisfying*

$$\sum_{k \in Z} |G_k(t)| \leq B, \quad a.e.$$

*but for some function $f \in L^2(R)$, the series*

$$\sum_{k \in Z} f(\cdot - k/b)G_k(\cdot)$$

*diverges in $L^2(R)$.*

*Proof.* We take $a = 1$, and $1/b$ irrational and let

$$f(t) = t^{-1/3}\chi_{(0,1/2]}(t).$$

We will construct $G_k$ of the form

$$G_k(t) = \sum_{n=-\infty}^{\infty} \chi_{I_k}(t+n),$$

where the $(I_k)$ are pairwise disjoint intervals contained in $(0, 1]$. We choose numbers $\epsilon_n$ with $0 < \epsilon_n < 1/2$, for $n = 1, 2, \cdots$ with

(6.3) $$\sum_n \epsilon_n = 1 \text{ and } \sum_n \epsilon_n^{1/3} = \infty.$$

We create intervals

$$J_n = (\sum_{k=1}^{n-1} \epsilon_k, \sum_{k=1}^{n} \epsilon_k] =: (a_n, b_n], \quad n = 1, 2, \cdots,$$

so that the $J_n$'s are pairwise disjoint with union $(0, 1]$. We define a sequence of integers $(k_n)_{n=0}^{\infty}$ satisfying $k_1 < k_2 < k_3 < \cdots$ in the following fashion. Let $k_0 = 1$, and when $n = 1, 2, \cdots$, we let $k_n$ be the smallest integer $k > k_{n-1}$ such that

$$\frac{k}{b} \pmod{Z} \in (a_n, \frac{a_n + b_n}{2}].$$

Next we define for $n = 1, 2, \cdots$,

$$c_n = \frac{k_n}{b} \pmod{Z}, \text{ and } I_{k_n} = J_n,$$

and let $I_k = \phi$, for all other $k$. Now, for all $n = 1, 2, \cdots$ we have

(6.4) $\|f(t - k_n/b)G_{k_n}(t)\|_{L^2(R)}^2 = \int_{c_n}^{b_n} f^2(t - c_n)dt \geq \int_0^{\frac{\epsilon}{2}} t^{-2/3}dt = 3(\frac{1}{2})^{1/3}\epsilon_n^{1/3}.$

Also, it follows that

$$\sum_k |G_k(t)| \leq 1, \quad \text{a.e.}$$

Moreover, the terms $f(t - k/b)G_k(t)$ in the series

$$\sum_k f(t - k/b)G_k(t)$$

have disjoint supports. Whence,

(6.5) $$\|\sum_{k \in Z} f(t - k/b)G_k(t)\|_{L^2(R)}^2 = \sum_{k \in Z} \|f(t - k/b)G_k(t)\|_{L^2(R)}^2.$$

It follows from (6.4) and (6.5) that the Walnut series for f diverges.

Now we consider the rational case (i.e. ab is rational) and discover that the convergence of the Walnut series becomes equivalent to the CC-condition. To simplify the proof, we first formulate and prove a lemma.

**Lemma 6.7.** *Let $(S, \Lambda, \mu)$ be a $\sigma$-finite measure space. Let $f_n \in L^\infty(S)$, for all $n \in Z$, and assume that the series $\sum_{n \in Z} h \cdot f_n$ converges unconditionally for every $h \in L^2(S)$. Then there is an $M > 0$ so that for all choices $b_n = 0, \pm 1$ and all $\theta_n \in R$ we have*
$$\|\sum_n b_n e^{i\theta_n} f_n\|_\infty \leq M.$$

*Proof.* For each choice of $b = (b_n)$ and $\theta = (\theta_n)$ and all $I \subset Z$ with $|I| < \infty$, define $T_{b,\theta,I} : L^2(S) \to L^2(S)$ by:
$$T_{b,\theta,I} h = \sum_{n \in I} h \cdot b_n e^{i\theta_n} f_n = h \cdot \sum_{n \in I} b_n e^{i\theta_n} f_n.$$

These operators are clearly bounded. By Proposition 2.6, our assumption in the lemma is just that the family of operators $(T_{b,\theta,I})$ are pointwise bounded. So by the Uniform Boundedness Principle, this family of operators has uniformly bounded norms. But, their norms are given by:
$$\|T_{b,\theta,I}\| = \|\sum_{n \in I} b_n e^{i\theta_n} f_n\|_\infty.$$

**Theorem 6.8.** *Let $ab$ be rational and assume that $(G_k)$ are periodic functions of period $a$. Then the following are equivalent:*

*(1) We have*
$$\sum_{k \in Z} |G_k(t)| \leq B, \quad a.e.$$

*(2) For every function $f \in L^2(R)$, the Walnut series*

(6.6)
$$\sum_{k \in Z} f(\cdot - k/b) G_k(\cdot),$$

*converges unconditionally in $L^2(R)$.*

*Proof.* (1) $\Rightarrow$ (2): Let $ab = p/q$ with $0 < p < q \in Z$ and gcd(p,q) $= 1$. After a change of scale, we may assume that $a = 1$. We re-group the Walnut series
$$\sum_k f(t - k/b) G_k(t) = \sum_k f(t - k\frac{q}{p}) G_k(t) = \sum_{j=0}^{p-1} \sum_\ell f(t - j\frac{q}{p} - \ell q) G_{j+\ell p}(t).$$

Each of the terms of the right-hand series above over $j$ has the form
$$\sum_\ell h(t - \ell q) H_\ell(t), \quad h \in L^2(R),$$

where the $H_\ell(t)$ are 1-periodic in $t$, and

$$\sum_\ell |H_\ell(t)| \leq B.$$

Now,

$$\int_{-\infty}^\infty \left[\sum_\ell |h(t-\ell q)H_\ell(t)|\right]^2 dt = \int_0^1 \sum_{n=-\infty}^\infty \left[\sum_\ell |h(t-\ell q+n)||H_\ell(t)|\right]^2 dt.$$

Take any $t \in [0,1]$ with

$$\sum_\ell |H_\ell(t)| \leq B \quad \text{and} \quad \sum_n |h(t-n)|^2 < \infty.$$

Setting

$$a_\ell = H_\ell(t) \quad \text{and} \quad z_\ell = (h(t-\ell q + n))_{n\in Z} \in \ell^2,$$

and using the triangle inequality, we have

(6.7) $$\|\sum_\ell a_\ell z_\ell\|_{\ell^2} \leq \sum_\ell |a_\ell| \|z_\ell\|_{\ell^2}.$$

Now, applying (6.7) pointwise, we have

$$\sum_{n=-\infty}^\infty \left(\sum_\ell |h(t-\ell q+n)||H_\ell(t)|\right)^2$$

$$\leq \left[\sum_\ell |H_\ell(t)| \left(\sum_{n=-\infty}^\infty |h(t-\ell q+n)|^2\right)^{1/2}\right]^2 \leq B^2 \sum_{n=-\infty}^\infty |h(t-n)|^2.$$

We can now conclude that

$$\int_{-\infty}^\infty \left(\sum_\ell |f(t-\ell q)||H_\ell(t)|\right)^2 \leq B^2 \int_0^1 \sum_n |h(t-n)|^2 dt = B^2 \|h\|^2.$$

It follows that the Walnut series is unconditionally convergent in $L^2(R)$-norm.

(2) $\Rightarrow$ (1): For this direction we use the Zak transform. So for $f \in L^2(R)$ and a.e. $t, \nu$, let

$$(Z_a f)(t,\nu) = a^{1/2} \sum_{\ell=1}^\infty f(a(t-\ell))e^{2\pi i \ell \nu}.$$

We compute,

$$(Z_a(G_k f(\cdot - k/b)))(t,\nu) = a^{1/2} \sum_{\ell=-\infty}^{\infty} G_k(a(t-\ell))f(a(t-\ell) - k/q)e^{-2\pi i \ell \nu}$$

$$= a^{1/2} G_k(at) \sum_{\ell=-\infty}^{\infty} f(a((t - k\tfrac{q}{p}) - \ell))e^{-2\pi i \ell \nu} = G_k(at)(Z_a f)(t - k\tfrac{q}{p}, \nu).$$

If we take $k = np + s$ with $n \in Z$ and $s = 0, 1, \cdots, p-1$, we can write the series in (6.6) in the Zak domain as

$$(6.8) \quad \sum_{k \in Z} (Z_a f)(t - k\tfrac{q}{p}, \nu) G_k(at) = \sum_{s=0}^{p-1} \sum_{n=-\infty}^{\infty} (Z_a f)(t - nq - s\tfrac{q}{p}, \nu) G_{np+s}(at)$$

$$= \sum_{s=0}^{p-1} (Z_a f)(t - s\tfrac{q}{p}, \nu) \sum_{n=-\infty}^{\infty} G_{np+s}(at) e^{-2\pi i n q \nu}.$$

Above we used the quasi-periodicity of the Zak transform in t (see (3.3)). By our assumptions and the unitarity of the Zak transform, we have that the right-hand side of (6.8) converges unconditionally for all $f \in L^2(R)$. Now let

$$(t, \nu) \in [0, 1) \times [0, 1) = \cup_{r=0}^{p-1} \left[\frac{r}{p}, \frac{r+1}{p}\right) \times [0, 1) =: \cup_{r=0}^{p-1} S_r.$$

Also, for $s, r = 0, 1, \cdots, p - 1$, let $U_{s,r}$ consist of all $f \in L^2(R)$ such that the restriction of $(Z_a f)(t - s\tfrac{q}{p}, \nu)$ to $[0, 1) \times [0, 1)$ is entirely supported by $S_r$. We note that for $s, r = 0, 1, \cdots, p - 1$, we have

$$\{(Z_a f)(\cdot - s\tfrac{q}{p}, \cdot) : f \in U_{s,r}\} \equiv L^2(S_r).$$

Thus, unconditional convergence of the right-hand side of (6.8) for all $f \in L^2(R)$, implies for all $h \in L^2(S_r)$ and all $s, r = 0, 1, \cdots, p-1$ the unconditional convergence of

$$h(t, \nu) \sum_{n=-\infty}^{\infty} G_{np+s}(at) e^{-2\pi i n q \nu}, \quad (t, \nu) \in S_r.$$

For $s, r = 0, 1, \ldots, p - 1$, by Lemma 6.7 we can find constants $M_{s,r}$ such that

$$(6.9) \quad \|\sum_{n \in Z} b_n e^{i\theta_n} G_{np+s}(at) e^{-2\pi i n q \nu}\|_{\infty, S_r} \leq M_{s,r},$$

for all choices $b_n = 0, 1$ and all $\theta_n \in R$. Let $t_0 \in [\frac{r}{p}, \frac{r+1}{p})$ be such that $at_0$ is a Lebesgue point of all $(G_{np+s})$. The set of these $t_0$ has full measure. We claim that

(6.10) $$\sum_{n \in Z} |G_{np+s}(at_0)| \leq M_{s,r}.$$

If inequality (6.10) does not hold, then we can find a finite set $I$ of n's such that

$$\sum_{n \in I} |G_{np+s}(at_0)| \geq M_{s,r}.$$

Take any $\nu_0 \in [0, 1)$ and choose $(\theta_n)$ with

$$e^{i\theta_n} e^{-2\pi i n q \nu_0} G_{np+s}(at_0) = |G_{np+s}(at_0)|.$$

Finally let

$$b_n = \begin{cases} 1 & : n \in I \\ 0 & : n \notin I. \end{cases}$$

Now,

$$|\sum_{n \in Z} b_n e^{i\theta_n} e^{-2\pi i n q \nu} G_{np+s}(at)| \geq M_{s,r},$$

in a $(t, \nu)$-set of positive measure, which contradicts (6.9). It now follows that

$$\operatorname{ess\,sup}_t \sum_{k=-\infty}^{\infty} |G_k(at)| \leq \max_{r=0,1,\cdots,p-1} \sum_{s=0}^{p-1} M_{s,r} =: M.$$

This completes the proof of Theorem 6.8.

We have seen that it takes special conditions to guarantee that the Walnut series converges unconditionally for all $f \in L^2(R)$. Now we will consider the more delicate question: If the conditions for unconditional convergence of the Walnut series for all $f \in L^2(R)$ are not satisfied can we at least find a "large class of functions" for which we do have unconditional convergence of the Walnut series? We will show that for $g \in \text{PF}$, the answer is yes. Furthermore, the "large class of functions" we exhibit will be independent of the choice of $g$.

**Notation 6.9.** *Let $f$ be any function on $R$, and $a, b \in R$. For any $k, \ell \in Z$ we define for all $t \in R$,*

(6.11) $$A_f(t) = \sup_{\ell \in Z} \sum_{k \in Z} |f(t - \ell a - k/b)|,$$

(6.12) $$B_f(t) = \sup_{k \in Z} \sum_{\ell \in Z} |f(t - \ell a - k/b)|.$$

**Theorem 6.10.** *Assume that $g \in PF$, $f \in L^2(R)$, and $A_f, B_f \in L^2([0, a))$. Then the series*

$$\sum_{k \in Z} |f(t - k/b)||G_k(t)|$$

*converges in norm in $L^2(R)$.*

*Proof.* We compute

$$\int_{-\infty}^{\infty} \left( \sum_{k \in Z} |f(t - k/b)||G_k(t)| \right)^2 dt = \int_0^a \sum_{\ell=-\infty}^{\infty} \left( \sum_{k \in Z} |f(t - \ell a - k/b)||G_k(t)| \right)^2 dt.$$

Take any $t \in [0, \infty)$ such that $A_f(t) < \infty$ and $B_f(t) < \infty$, and let $B$ be as in Proposition 2.3. By the Schur Test, stated before Proposition 4.6, the matrix

$$(|f(t - \ell a - k/b)|)_{\ell, k \in Z}$$

defines a bounded linear operator on $\ell^2(Z)$ with operator norm $\leq \max(A_f(t), B_f(t))$. Hence,

(6.13)
$$\sum_{\ell=-\infty}^{\infty} \left( \sum_{k \in Z} |f(t - \ell a - k/b)||G_k(t)| \right)^2$$

$$\leq \max^2(A_f(t), B_f(t)) \sum_k |G_k(t)|^2 \leq B \max^2(A_f(t), B_f(t)).$$

Since (6.13) holds a.e. $t \in R$, we have

$$\int_{-\infty}^{\infty} \left( \sum_{k \in Z} |f(t - k/b)||G_k(t)| \right)^2 dt \leq B \int_0^a \max^2(A_f(t), B_f(t)) dt < \infty,$$

by the assumptions in the theorem.

Note that the conditions in Theorem 6.10 are reasonably weak, and that the condition on the function $f$ is independent of $g$. In particular, we thus get a generalization of Proposition 2.4 in which now also non-compactly supported f's, such as $f \in W(L^\infty, L^1)$, are allowed.

We end this section with a discussion concerning generalizations of Propositions 5.3, 5.6, and 6.2 to the case where $ab$ is rational. Let $(g, a, b)$ be a WH-system with a finite upper frame bound and assume $ab = p/q$ with $p, q \in Z$ and gcd(p,q) = 1. We define for $m = 0, 1, \cdots, p - 1$,

(6.14) $F_{km}(t,\nu) =$

$$\frac{1}{p}e^{-2\pi ik(\nu+m/p)}\sum_{j=0}^{q-1}\sum_{r=0}^{p-1}e^{-2\pi irm/p}\cdot\sum_{d=-\infty}^{\infty}g\left(\frac{t-r-\frac{jp}{q}-dp}{b}\right)\overline{g\left(\frac{t-r-\frac{jp}{q}-dp}{b}\right)}.$$

We consider the Hilbert space of column vectors of the form

$$\begin{bmatrix} f_0 \\ f_1 \\ \vdots \\ f_{p-1} \end{bmatrix},$$

where $f_r \in L^2([0,1) \times [0,1/p))$ for $r = 0, 1, \cdots, p-1$, and norm given by

$$\left\|\begin{bmatrix} f_0 \\ f_1 \\ \vdots \\ f_{p-1} \end{bmatrix}\right\|^2 = \sum_{r=0}^{p-1}\int_0^1\int_0^{1/p}|f_r(t,\nu)|^2 dt\, d\nu.$$

Now, $L^2(R)$, $\|\cdot\|_{L^2}$ and the above Hilbert space can be identified with one-another according to

$$f \in L^2(R) \to \begin{bmatrix} (Z_{1/b}f)(t,\nu+\frac{0}{p}) \\ \vdots \\ (Z_{1/b}f)(t,\nu+\frac{p-1}{p}) \end{bmatrix},$$

for a.e. $t, \nu \in [0,1) \times [0,1/p)$. That is, we set for $r = 0, 1, \cdots, p-1$,

$$f_r(t,\nu) = (Z_{1/b}f)(t,\nu+\frac{r}{p}), \quad \text{a.e.} \quad t,\nu \in [0,1) \times [0,1/p).$$

Now we will compute

(6.15) $$Z_{1/b}[G_k(\cdot)f(\cdot - k/b)](t,\nu).$$

We have $a = \frac{p}{q}\frac{1}{b}$. Thus we can replace $G_k(t)$ by taking $n = dq+j$, $j = 1, 2, \cdots, q-1$, $d \in Z$, so that

(6.16) $$G_k(t) = \sum_{j=0}^{q-1}\sum_{d=-\infty}^{\infty}g(t - \frac{jp}{qb} - \frac{dp}{b})\overline{g(t - \frac{jp}{qb} - \frac{dp}{b} - \frac{k}{b})}.$$

Therefore,

$$Z[G_k(\cdot)f(\cdot - k/b)(t,\nu) = b^{-1/2} \sum_{\ell=-\infty}^{\infty} G_k(\frac{t-\ell}{b})f(\frac{t-\ell-k}{b})e^{2\pi i\ell\nu}$$

$$= b^{-1/2} \sum_{\ell=-\infty}^{\infty} \sum_{j=0}^{q-1} \sum_{d=-\infty}^{\infty} g(\frac{t-\ell-\frac{jp}{q}-dp}{b})\overline{g(\frac{t-\ell-\frac{jp}{q}-dp}{b})}f(\frac{t-\ell-k}{b})e^{2\pi i\ell\nu}.$$

Next we let $\ell = np + r$, $n \in Z$, $r = 0, 1, \cdots, p-1$, so that

$$Z[G_k(\cdot)g(\cdot - k/b)](t,\nu) =$$

$$b^{-1/2} \sum_{j=0}^{q-1}\sum_{r=0}^{p-1} \left[ \sum_{d=-\infty}^{\infty} g(\frac{t-r-\frac{jp}{q}-(d+n)p}{b})\overline{g(\frac{t-r-k-\frac{jp}{q}-(d+n)p}{b})} \right] \cdot$$

$$(\sum_{n=-\infty}^{\infty} f(\frac{t-r-k-np}{b})e^{2\pi i(np+r)\nu}).$$

Next we write

(6.17) $$\sum_{n=-\infty}^{\infty} f(\frac{t-r-k-np}{b})e^{2\pi i(np+r)\nu} =$$

$$\frac{1}{p} \sum_{m=0}^{p-1} \sum_{\ell=\infty}^{\infty} f(\frac{t-k-\ell}{b})e^{2\pi i\ell\nu + 2\pi i(\ell-r)m/p},$$

where in (6.17) we have used the formula

$$\frac{1}{p} \sum_{m=0}^{p-1} e^{2\pi ijm/p} = \begin{cases} 1 & \text{for j a multiple of p} \\ 0 & \text{otherwise.} \end{cases}$$

Now we get

$$Z[G_k(\cdot)f(\cdot - k/b)](t,\nu) =$$

$$\frac{1}{p} \sum_{j=0}^{q-1}\sum_{r=0}^{p-1}\sum_{m=0}^{p-1} \left[ g(\frac{t-r-\frac{jp}{q}-dp}{b})\overline{\frac{t-r-k-\frac{jp}{q}-dp}{b}}) \right] \cdot$$

$$e^{-2\pi ik(\nu+m/p)-2\pi irm/p}(Z_{1/b}f)(t,\nu+m/p)$$

$$= \sum_{m=0}^{p-1} F_{km}(t,\nu)(Z_{1/b}f)(t,\nu+m/p),$$

where

$$F_{km}(t,\nu) = \frac{1}{p}e^{-2\pi ik(\nu+m/p)}\sum_{j=0}^{q-1}\sum_{r=0}^{p-1}e^{-2\pi irm/p}\cdot\sum_{d=-\infty}^{\infty}g(\frac{t-r-\frac{jp}{q}-dp}{b})\overline{g(\frac{t-r-\frac{jp}{q}-dp}{b})}.$$

It follows that

(6.18)
$$Z_{1/b}[\sum_{k=-L}^{K}G_k(\cdot)f(\cdot-k/b)](t,\nu) =$$

$$\sum_{m=0}^{p-1}(\sum_{k=-L}^{K}F_{km}(t,\nu))(Z_{1/b}f)(t,\nu+\frac{m}{p}), \quad \text{a.e. } t,\nu \in R.$$

Assume that the Walnut series converges strongly for all $f \in L^2(R)$. By taking f's such that $(Z_{1/b}f)(t,\nu+m/p)$ has support in $[0,1) \times [0,1)$ entirely contained in one strip $[0,1) \times [\frac{r}{p}, \frac{r+1}{p})$ with $r = 0, 1, \cdots, p-1$, we see from (6.18) that for all $m = 0, 1, \cdots, p-1$

(6.19)
$$\sum_{k=-L}^{K}F_{km}(t,\nu)$$

has a uniformly bounded $L^\infty([0,1) \times [0,1)$-norm. conversely, assume that (6.19) has a uniformly bounded $L^\infty([0,1) \times [0,1)$-norm for $m = 0, 1, \cdots, p-1$. Then it follows from Theorem 5.5 and (6.18) that the Walnut series converges strongly for every $f \in L^2(R)$. This generalizes Proposition 5.6, and the Propositions r.3 and 6.2 can be generalized similarly. Thus we get

**Proposition 6.11.** *Let* $ab = p/q$ *with* $gcd(p,q) = 1$ *and let* $g \in PF$.
*(I) The following are equivalent:*
*(1) The Walnut series converges symmetrically for every* $f \in L^2(R)$.
*(2) There is a* $B > 0$ *so that*

$$sup_K ess\ sup_{t,\nu}|\sum_{k=-K}^{K}F_{km}(t,\nu)| \leq B,$$

*for all* $m = 0, 1, \cdots, p-1$.
*(II) The following are equivalent:*
*(3) The Walnut series converges strongly (i.e. in norm) for all* $f \in L^2(R)$.

(4) There is a $B > 0$ so that

$$sup_{K,L}\, ess\, sup_{t,\nu} |\sum_{k=-K}^{L} F_{km}(t,\nu)| \leq B,$$

for all $m = 0, 1, \cdots, p - 1$.

(III) The following are equivalent:

(5) The Walnut series converges unconditionally for every $f \in L^2(R)$.

(6) There is a $B > 0$ so that

$$sup_M\, ess\, sup_{t,\nu} |\sum_{k \in M} F_{km}(t,\nu)| \leq B,$$

for all $m = 0, 1, \cdots, p - 1$, and the first sup runs over all finite sets $M \subset Z$.

Finally, note that when $a = b = p = q = 1$, we only have to consider the case m=0 and (6.14) becomes,

$$F_{k0}(t,\nu) = e^{-2\pi i k\nu} \sum_{d=-\infty}^{\infty} g(t - d)\overline{g(t - k - d)} = e^{-2\pi i k\nu} G_k(t).$$

That is, the sums $\sum_k e^{-2\pi i k\nu} G_k(t)$ which we consider in Propositions 5.3, 5.6, and 6.2 are just the sums $\sum_k F_{km}(t,\nu)$.

## 7. Extending the Frame Operator

In this section we will classify when the frame operator extends to a bounded operator on other classes of spaces. We will see that the CC-condition, although much weaker than previously used conditions which imply convergence of the frame operator, is also a very strong assumption on the function $g$.

**Theorem 7.1.** *If $ab \leq 1$ and $g \in L^2(R)$, the following are equivalent:*

(1) There is a constant $B > 0$ so that

$$\sum_{k \in Z} |G_t(x)| \leq B, \quad a.e. \quad x \in R.$$

(2) The frame operator $Sf = \sum_{n,m \in Z} <f, E_{mb}T_{na}g> E_{mb}T_{na}g$ extends to a bounded linear operator from $L^p(R)$ to $L^p(R)$ for every $1 \leq p \leq \infty$. (Here, by $L^\infty$, we really mean $L_0^\infty$, the closure of the compactly supported functions in $L^\infty$).

*Moreover, if ab is rational and g satisfies the uniform CC-condition, then S extends to an isomorphism of $L^p(R)$ to $L^p(R)$ for every $1 \leq p \leq \infty$.*

*Proof.* (1) $\Rightarrow$ (2): We will show that $S$ is a bounded linear operator mapping a dense subset of $L^1$ into $L^1$ and a dense subset of $L_0^\infty$ into $L_0^\infty$. By the Riesz-Thorin interpolation theorem [27], p. 95, $S$ then extends to a bounded linear operator from $L^p$ to $L^p$ for all $1 \leq p < \infty$.

**Case I: The $L^1$-Case.**

If $f \in L^1$ is bounded and compactly supported then by Proposition 2.4 we have:

$$\|Sf\|_{L^1} = \|Lf\|_{L^1} = b^{-1} \int_R |\sum_k (T_{k/b}f)G_k|$$

$$\leq b^{-1} \int_R \sum_k |(T_{k/b}f)G_k| = b^{-1} \sum_k \int_R |(T_{k/b}f)G_k|$$

$$= b^{-1} \sum_k \int_R |f||\overline{G_{-k}}| \leq b^{-1} \int_R |f| \sum_k |G_k| \leq b^{-1} K \int_R |f| = b^{-1} B \|f\|_{L^1}.$$

The above makes $S$ a bounded linear operator from a dense subspace of $L^1$ into itself, and hence it uniquely extends to a bounded linear operator on $L^1$.

**Case II: $L_0^\infty$-Case.**

For any $f \in L^\infty$ which is compactly supported we have:

$$\|Sf\|_{L^\infty} = \|Lf\|_{L^\infty} = \|\sum_k (T_{k/b}f)G_k\|_{L^\infty} = \text{ess sup} |\sum_k (T_{k/b}f)G_k|$$

$$\leq \text{ess sup} \sum_k |T_{k/b}f||G_k| \leq \text{ess sup} \sum_k \|f\|_{L^\infty}|G_k|$$

$$\leq \|f\|_{L^\infty} \text{ess sup} \sum_k |G_k| \leq B\|f\|_{L^\infty}.$$

Again, this makes $S$ a bounded linear operator on a dense subset of $L_0^\infty$ and hence it uniquely extends to a bounded linear operator on $L_0^\infty$.

(2) $\Rightarrow$ (1): By our assumption, $S$ is a bounded linear operator on $L_0^\infty$. Fix n, let $I = [0, a]$ and choose functions $(f_j)_{j=-n}^n$ satisfying:
(1) $|f_j| = \chi_{[j/b,(j/b)+a]}$
(2) $(T_{-j/b}f_j)G_j = \chi_{[0,a]}|G_j|$.

Let $f = \sum_{j=-n}^{n} f_j$. Then $\|f\|_{L^\infty} = 1$, here we use the assumption that $ab \leq 1$, and so

$$\|S\| \geq \|SF\|_{L^\infty} \geq \|\chi_{[0,a]} Sf\|_{L^\infty} = \|\sum_{k=-n}^{n} (T_{-k/b} f_k) G_k\|_{L^\infty}$$

$$= \|\sum_{k=-n}^{n} |G_k|\|_{L^\infty} = \text{ess sup} \sum_{k=-n}^{n} |G_k|.$$

We dropped the restriction to $[0, a]$ in the last step since all the $G_k$ are periodic of period a. Since n was arbitrary, we have (1).

For the moreover part of the theorem, we apply Theorem 4.14 to conclude that $S^{-1}g$ (which has $S^{-1}$ as its frame operator) satisfies the CC-condition and apply the first part of the theorem to $S^{-1}$.

Next we will classify when the frame operator extends to a bounded linear operator on the Wiener amalgam space $W(L^\infty, \ell^1)$. We will need the properties of this space contained in Lemma 2.5. A good reference for this topic is Chapter 3 in [10], due to Feichtinger and Zimmermann.

**Theorem 7.2.** *If $ab \leq 1$ and $g \in PF$, the following are equivalent:*

*(1) The frame operator is a bounded linear operator from $W(L^\infty, \ell^1)$ to $W(L^\infty, \ell^1)$.*

*(2) We have*

$$\sum_{k \in Z} \|G_k\|_\infty = \sum_{k \in Z} \text{ess sup} |G_k(x)| = B < \infty.$$

*Proof.* (1) $\Rightarrow$ (2): By Lemma 2.5 and the fact that $a \leq 1/b$ and $G_k$ is periodic of period a, we have for all $f = \chi_{[0,a]}$

$$2\|Sf\|_{W,a} \geq \|Sf\|_{W,1/b} = \sum_k \|\chi_{k/b,(k/b)+a]} G_k\|_\infty = \sum_k \|G_k\|_\infty.$$

Hence,

$$\sum_{k \in Z} \|G_k\|_\infty \leq 2\|S\| \|\chi_{[0,a]}\|_\infty = 2\|S\|.$$

(2) $\Rightarrow$ (1): If we have (2) then the CC-condition is satisfied and so the Walnut series for Sf converges unconditionally by Corollary 6.5. Hence, for any $f \in W(L^\infty, \ell^1)$ if we let

$$f_j = \chi_{[j/b,(j+1)/b]} f$$

then

$$\|f\|_{W,1/b} = \sum_{j \in Z} \|f_j\|_\infty.$$

It follows that
$$b\|Sf\|_{W,1/b} = \|\sum_k (T_{k/b}f)G_k\|_{W,1/b} =$$
$$\sum_{j\in Z}\|\sum_k \chi_{[j/b,(j+1)/b]}(T_{k/b}f)G_k\|_\infty \le \sum_{j\in Z}\sum_{k\in Z}\|\chi_{[j/b,(j+1)/b]}T_{k/b}f\|_\infty \|G_k\|_\infty =$$
$$\sum_{j\in Z}\sum_{k\in Z}\|f_{j+k}\|_\infty \|G_k\|_\infty = \sum_{k\in Z}\|G_k\|_\infty \sum_{j\in Z}\|f_{j+k}\|_\infty = \|f\|_{W,1/b}\sum_{k\in Z}\|G_k\|_\infty.$$

If we choose a natural number $m$ so that $1/b \le ma$ we now have by Lemma 2.5,
$$\|Sf\|_{W,a} \le 2m\|Sf\|_{W,1/b} \le 2mb^{-1}\sum_k \|G_k\|_\infty \|f\|_{W,1/b} \le 4mb^{-1}\sum_k \|G_k\|_\infty \|f\|_{W,a}.$$

That is, $S$ is a bounded linear operator on $W(L^\infty, \ell^1)$.

The assumption that $g \in W(L^\infty, \ell^1)$ is so strong that we can easily obtain results for $S^{-1}$ in this case.

**Corollary 7.3.** *If $g \in W(L^\infty, \ell^1)$, then $(E_{mb}T_{na}g)$ has a finite upper frame bound for all $ab \le 1$ and the frame operator for this frame is a bounded linear operator on $W(L^\infty, \ell^1)$. Moreover, in this case for each $a \in R$ there is a $b_0 \in R$ so that for all $0 < b \le b_0$, if $G_{0,a,b}$ is bounded below then $(g, a, b)$ is a frame whose frame operator is an isomorphism on $W(L^\infty, \ell^1)$.*

*Proof.* Given any $f, h \in W(L^\infty, \ell^1)$ we have
$$\sum_k \|\sum_n |T_{na}f||T_{na+k/b}h|\|_\infty = \sum_k \|\sum_n |T_{na}f||T_{na+k/b}h|\chi_{[0,a]}\|_\infty$$
$$\le \sum_k \sum_n \|(T_{na}f)\chi_{[0,a]}\|_\infty \|(T_{ka+b/k}h)\chi_{[0,a]}\|_\infty$$
$$\le \sum_n \|(T_{na}f)\chi_{[0,a]}\|_\infty \sum_k \|(T_{k/b}(T_{na}h))\chi_{[0,1/b]}\|_\infty$$
$$\le \sum_n \|(T_{na}f))\chi_{[0,a]}\|_\infty \|T_{na}h\|_{W,1/b}$$
$$\le 2\|h\|_{W,1/b}\sum_n \|(T_{na}f)\chi_{[0,a]}\|_\infty \le 4\|h\|_{W,a}\|f\|_{W,a}.$$

If $g \in W(L^\infty, \ell^1)$, we can let $f = h = g$ above and get that
$$\sum_k \|G_k\|_\infty \le 4\|g\|_{W,a}^2.$$

By Theorem 7.2, we have the first part of the Corollary. The second part comes from the Walnut representation for $S^{-1}$ and similar calculations. (See Heil and Walnut [14], Theorem 4.2.2).

## 8. Classifying WH-Systems with the Same Frame Operator

In this section we will classify those functions which induce equal frame operators. As we will see, this is a natural generalization of a result of Ron and Shen [23], Corollary 2.19 which we state below. See also [15], Subsection 1.3, [3], Theorem 4.2, or [11] for a generalization.

**Theorem 8.1.** *Let $a, b \in R$ and $g \in PF$. The following are equivalent:*
(1) $(E_{mb}T_{na}g)$ *is a normalized tight Weyl-Heisenberg frame for $L^2(R)$.*
(2) *We have:*
   (a) $G_g(t) = \sum_{n \in Z} |g(t - na)|^2 = b$ *a.e.,*
   (b) $G_{g,k}(t) = \sum_{n \in Z} g(t - na)\overline{g(t - na - k/b)} = 0$ *a.e. for all $k \neq 0$.*

One interpretation of Theorem 8.1 is that two normalized tight Weyl-Heisenberg frames $(E_{mb}T_{na}g)$ and $(E_{mc}T_{nd}h)$ have the same frame operator - in this case the identity operator - if and only if $b^{-1}G_g = d^{-1}G_h$ a.e. and for all $0 \neq k \in Z$ we have that $0 = G_{g,k} = G_{h,k}$ a.e. Our next result is a natural extension of this to the case of arbitrary frames.

**Theorem 8.2.** *Let $(g, a, b)$ and $(h, c, d)$ be preframe WH-systems. Then the following are equivalent:*
(1) *Their frame operators are equal, i.e. $S_g = S_h$.*
(2) *One of the following holds:*
   (i) $\frac{d}{b}$ *is not rational. Then $b^{-1}G_{g,0} = d^{-1}G_{h,0}$ a.e. and for all $0 \neq k \in Z$ we have*
$$G_{g,k} = G_{h,k} = 0 \quad a.e.$$
   (ii) $\frac{d}{b} = \frac{p}{q}$ *is rational where $p, q$ are natural numbers. Then we have for all $k \in Z$*
$$b^{-1}G_{g,qk} = d^{-1}G_{h,pk} \quad a.e.$$
*and for all other integers $m \neq qk$ and $\ell \neq pk$ we have*
$$G_{g,m} = G_{h,\ell} = 0 \quad a.e.$$

*Moreover, in all the above cases, either there are natural numbers $r < s$ so that $a = \frac{r}{s}c$ (or $c = \frac{r}{s}a$) and all $G_{g,k}$, $G_{h,k}$ are periodic of period $1/s$, or $\frac{a}{c}$ is irrational and all $G_{g,k}$, $G_{h,k}$ are constant a.e.*

*Proof.* (1) $\Rightarrow$ (2): We assume (1) is true and check the two cases of (2) separately.

**Case I.** *We assume that $\frac{d}{b}$ is irrational.*

We are assuming that $S_g = S_h$. We may also assume that $b < d$. Fix $k \in Z$ and for any $0 \neq e \in R$, let $E_e = \{k/e : k \in Z\}$. Fix $I \subset R$ with $|I| < 1/d$ and $I \cap E_b = \phi = I \cap E_d$. Then for all bounded $f \in L^2(I)$, since $S_g f = S_h f$, invoking the Walnut representation we have

$$0 = \chi_I[S_g f - S_h f] = f \, \chi_I[b^{-1}G_{g,0} - d^{-1}G_{h,0}].$$

It follows easily from here that

$$b^{-1}G_{g,0} = d^{-1}G_{h,0} \quad \text{a.e.}$$

Next, fix $0 \neq k \in Z$ and let

$$\epsilon = \min \{|\frac{-k}{d} - \frac{\ell}{b}| : \ell \in Z\}.$$

Our assumption that $\frac{d}{b}$ is irrational implies that $\epsilon > 0$. For any interval $I \subset [0, 1/d]$ with $|I| < \epsilon$ let

$$f = \chi_{(\frac{-k}{d}+I)}.$$

Then $\chi_{[0,1/d]}S_h f = \chi_I$ while our assumption that $|I| < \epsilon$ implies that $\chi_{[0,1/d]}S_g f$ has no support in the interval I. Since we are assuming that $S_g = S_h$, we conclude from

$$0 = \chi_{[0,1/d]}(S_g f - S_h f),$$

that

(8.1) $$0 = \chi_{[0,1/d]}S_g f.$$

If we expand $S_g f$ via its Walnut representation, (8.1) becomes

$$0 = \chi_I d^{-1} G_{h,k}.$$

Thus $G_{h,k} = 0$ a.e. Now, for any natural number $k > 0$, let

$$f = \chi_{[\frac{-k}{d}, \frac{-(k-1)}{d}]}.$$

Since $G_{h,k} = 0$ a.e., we have that

$$0 = \chi_{[0,1/b]}[S_g f - S_h(f)] = \chi_{[0,1/b]}G_{g,k}.$$

Hence, $G_{g,k} = 0$ a.e., for all natural numbers k. Similarly, $G_{g,k} = 0$ a.e. for all integers $k < 0$.

**Case II.** *We assume that $\frac{d}{b} = \frac{p}{q}$ with $\frac{p}{q}$ minimal and in lowest form.*

Exactly the same argument as in Case I shows that $G_{h,\ell} = 0$ a.e. for all $\ell \neq pk$, $k \in Z$. Again, the second half of the argument of Case I shows now that $G_{g,m} = 0$ a.e. for all $m \neq qk$, $k \in Z$. Now, fix $0 < k \in Z$ and let

$$f = \chi_{[\frac{-qk}{b}, \frac{q(k-1)}{b}]}.$$

Using the fact that $G_{h,\ell} = 0$ a.e. for all $\ell \neq pk$, $k \in Z$, and $G_{g,m} = 0$ a.e. for all $m \neq qk$, $k \in Z$ we see that

$$0 = \chi_{[0,1/b]}[S_g f - S_h f] = \chi_{[0,1/b]}[b^{-1} G_{g,qk} - d^{-1} G_{h,pk}] \text{ a.e.}$$

It follows that $b^{-1} G_{g,qk} = d^{-1} G_{h,pk}$ a.e. on $[0, 1/b]$. But, all $G_{g,k}$ are periodic of period $a < 1/b$ and all $G_{h,k}$ are periodic of period $c < 1/d < 1/b$. Therefore, these functions are equal a.e. on all of $R$.

$(2) \Rightarrow (1)$: Again we will check that each of the two cases implies $(1)$.

**Case I.** *We have that $\frac{d}{b}$ is not rational.*

For any bounded, compactly supported $f \in L^2(R)$ we have by (ii),

$$S_g(f) - S_h(f) = b^{-1} \sum_k T_{k/b}(f) G_{g,k} - d^{-1} \sum_k T_{k/d}(f) G_{h,k}$$

$$= b^{-1} f \, G_{g,0} - d^{-1} f \, G_{h,0} = 0.$$

Again we conclude that $S_g = S_h$ on $L^2(R)$.

**Case II.** *We have that $\frac{q}{b} = \frac{p}{d}$.*

For any bounded, compactly supported $f \in L^2(R)$ we have by (i),

$$S_g(f) - S_h(f) = b^{-1} \sum_{k \in Z} T_{k/b}(f) G_{g,k} - \sum_{k \in Z} T_{k/d}(f) G_{h,k}$$

$$= \sum_k T_{qk/b}(f) G_{g,qk} - \sum_k T_{pk/d}(f) G_{h,pk} = \sum_{k \in Z} T_{qk/b}(f)[G_{g,qk} - G_{h,pk}] = 0.$$

Since $S_g$ and $S_h$ are both bounded operators and are equal on a dense subset of $L^2(R)$, they are equal.

The moreover part of the theorem is a well-known result about functions which are periodic with two different periods a,c.

As a special case of the theorem when $a = c$ and $b = d$ we have exactly the obvious generalization of the tight frames case.

**Corollary 8.3.** *Let $(g, a, b)$ and $(h, a, b)$ be preframe WH-systems. Then the following are equivalent:*

*(1) $S_g = S_h$*

*(2) For all $k \in Z$ we have*

$$G_{g,k} = G_{h,k} \quad a.e.$$

DEPARTMENT OF MATHEMATICS, THE UNIVERSITY OF MISSOURI, COLUMBIA, MISSOURI 65211, USA; MATHEMATICAL INSTITUTE, BUILDING 303, TECHNICAL UNIVERSITY OF DENMARK, 2800 LYNGBY, DENMARK; PHILIPS RESEARCH LABORATORIES EINDHOVEN, 5656 AA EINDHOVEN, THE NETHERLANDS

*E-mail address*: pete@casazza.math.missouri.edu; olechr@mat.dtu.dk; janssena@natlab.research.philips.com